\documentclass[10pt,a4paper]{article}

\usepackage{amsmath,amsfonts,amssymb}

\newtheorem{defn}{Definition}[section]
\newtheorem{theo}[defn]{Theorem}
\newtheorem{lem}[defn]{Lemma}
\newtheorem{prop}[defn]{Proposition}
\newtheorem{cor}[defn]{Corollary}
\newtheorem{rem}[defn]{Remark}

\newenvironment{proof}{{\bf Proof }}{{\par\nobreak\vskip 0.1cm\nobreak\hfill$\Box$}}

\begin{document} 

\noindent
{\Large \bf A weak maximum principle and solvability for two-dimensional elliptic equations with $L^2$-drifts and singular zero-order coefficients}
\\ \\
\bigskip
\noindent
{\bf Haesung Lee}  \\
\noindent
{\bf Abstract.} We study homogeneous Dirichlet problems for second-order linear elliptic equations in divergence form on bounded open subsets $U$ of $\mathbb R^2$. The diffusion matrix is bounded, measurable, and uniformly elliptic, and the drift belongs to $L^2(U,\mathbb R^2)$. If the zero-order coefficient $c$ belongs to $L^1(U)$ and satisfies $c\ge0$, we prove a weak maximum principle and deduce uniqueness of $H^{1,p}$-weak solutions for every $p \in (2, \infty)$. No regularity of the boundary is required for this result. The proof uses a logarithmic test function and an estimate for circular means of Sobolev functions. If, in addition, $U$ is a bounded open set with Lipschitz boundary and $c\in L^s(U)$ for some $s>1$, then for every $q\in(2, \infty)$ there exists $\varepsilon\in(0,q-2)$ such that the Dirichlet problem admits a unique $H^{1,2+\varepsilon}$-weak solution for every right-hand side in $H^{-1,q}(U)$. For existence, we use estimates for the lower-order terms, the global $H^{1,p}$-estimate for the diffusion operator, and the method of continuity. Finally, on the unit disk $B_1\subset\mathbb R^2$, we construct $c\in L^1(B_1)$ with $c \geq 0$ and $f\in C_0^\infty(B_1)$ such that the Dirichlet problem for $-\Delta u+cu=f$ has a bounded weak solution in $H^{1,2}_0(B_1)$ but no weak solution in $H^{1,p}_0(B_1)$ for any $p\in(2,\infty)$.


\section{Introduction}\label{intro}

In this paper, we study the existence and uniqueness of $H^{1,p}$-weak solutions, with $p\in(2,\infty)$, to the Dirichlet problem
\begin{equation}\label{introproblem}
\left\{
\begin{aligned}
-\operatorname{div}(A\nabla u)+\langle\mathbf H,\nabla u\rangle+cu&=\ell
&&\text{in }U,\\
u&=0
&&\text{on }\partial U.
\end{aligned}
\right.
\end{equation}
Here, $U$ is a bounded open subset of $\mathbb R^2$, and $A=(a_{ij})_{1\le i,j\le2}$ is a possibly non-symmetric matrix of bounded measurable functions on $\mathbb R^2$. We assume that there exist constants $\lambda,M>0$ such that
\begin{equation}\label{elliptbdda}
\langle A(x)\xi,\xi\rangle\ge\lambda\|\xi\|^2,
\qquad
\max_{1\le i,j\le2}|a_{ij}(x)|\le M, \quad \text{ for a.e. $x\in\mathbb R^2$ and every $\xi\in\mathbb R^2$.}
\end{equation}
The lower-order coefficients satisfy
$$
\mathbf H\in L^2(U,\mathbb R^2),
\qquad
c\in L^1(U),
\qquad
c\ge0
\quad\text{a.e. in }U.
$$
For $q \in (1,\infty)$, we write
$$
q':=\frac{q}{q-1},
\qquad
H^{-1,q}(U):=\bigl(H^{1,q'}_0(U)\bigr)^*.
$$
For $p\in(2,\infty)$ and $\ell\in\mathcal D'(U)$, a function $u\in H^{1,p}_0(U)$ is called an $H^{1,p}$-weak solution of \eqref{introproblem} if
\begin{equation}\label{introweakformulation}
\int_U\langle A\nabla u,\nabla\varphi\rangle\,dx+\int_U\langle\mathbf H,\nabla u\rangle\varphi\,dx+\int_U cu\varphi\,dx=\langle\ell,\varphi\rangle \qquad\text{for every }\varphi\in C_0^\infty(U).
\end{equation}
Here, $\mathcal D'(U)$ denotes the space of distributions on $U$, and $\langle\ell,\varphi\rangle$ denotes the action of $\ell$ on $\varphi \in C_0^{\infty}(U)$. Since $p\in(2,\infty)$, the Sobolev embedding gives $u\in L^\infty(U)$, and hence $cu\in L^1(U)$. Moreover, since $U$ is bounded, $\nabla u\in L^2(U,\mathbb R^2)$, so H\"older's inequality yields $\langle\mathbf H,\nabla u\rangle\in L^1(U)$. Thus all terms in \eqref{introweakformulation} are well defined.\\
We summarize the uniqueness and existence results in the following theorem.
\begin{theo}[Main result]\label{thmmain}
Let $U\subset\mathbb R^2$ be a bounded open set, and let $A=(a_{ij})_{1\le i,j\le2}$ be a matrix of bounded measurable functions on $\mathbb R^2$ satisfying \eqref{elliptbdda}. Assume that
$$
\mathbf H\in L^2(U,\mathbb R^2),
\qquad
c\in L^1(U),
\qquad
c\ge0
\quad\text{a.e. in }U.
$$
Then the following assertions hold.
\begin{itemize}
\item[(i)] {\rm(Uniqueness).}
For every $p\in(2,\infty)$, if $u\in H^{1,p}_0(U)$ is an $H^{1,p}$-weak solution of \eqref{introproblem} with $\ell=0$, then $u=0$.
\item[(ii)] {\rm(Existence and uniqueness).}
Suppose, in addition, that $U$ is a bounded open set with Lipschitz boundary and $c\in L^s(U)$ for some $s\in(1,\infty)$. Then, for every $q\in(2,\infty)$, there exists $\varepsilon\in(0,q-2)$, depending only on $q$, $s$, the constants $\lambda,M$ in \eqref{elliptbdda}, and $U$, such that, with $p:=2+\varepsilon$, problem \eqref{introproblem} admits a unique $H^{1,p}$-weak solution $u\in H^{1,p}_0(U)$ for every $\ell\in H^{-1,q}(U)$. Moreover,
\begin{equation}\label{mainestimate}
\|u\|_{H^{1,p}(U)}\le C\|\ell\|_{H^{-1,q}(U)},
\end{equation}
where $C>0$ is a constant independent of $\ell$.
\end{itemize}
\end{theo}
Part {\rm(i)} requires no regularity of $\partial U$. By linearity, it yields uniqueness of $H^{1,p}$-weak solutions for any fixed right-hand side $\ell\in\mathcal D'(U)$ without asserting existence. In part {\rm(ii)}, the exponent $q\in(2, \infty)$ of the right-hand side may be chosen arbitrarily, while the solution exponent $p=2+\varepsilon$ is chosen in the admissible range specified in Theorem~\ref{thmexistence}. Since $2<p<q$ and $U$ is bounded, the continuous embedding $H^{-1,q}(U)\hookrightarrow H^{-1,p}(U)$ allows us to apply the isomorphism result of that theorem at the exponent $p$. Theorem~\ref{thmexistence} also provides estimates for solutions with right-hand sides of the form $f-\operatorname{div}\mathbf F$. We prove part {\rm(i)} in Section~\ref{uniquenesec} and part {\rm(ii)} in Section~\ref{existsec}.\\
The homogeneous Dirichlet problem for elliptic equations with measurable coefficients was studied extensively in \cite{Stampacchia65,Trudinger73}, where maximum principles, existence and uniqueness results, and regularity estimates were established under suitable assumptions on the lower-order coefficients. In the usual Lebesgue-space setting, the drift assumption $\mathbf H\in L^d(U,\mathbb R^d)$ allows the drift term to be treated in $H^{1,2}_0(U)$ when $d\ge3$, whereas $\mathbf H\in L^q(U,\mathbb R^2)$ with $q>2$ is typically assumed when $d=2$. This difference comes from the Sobolev embedding used to estimate the drift term. For $d\ge3$, the embedding
$$
H^{1,2}_0(U)\hookrightarrow L^{\frac{2d}{d-2}}(U)
$$
allows one to control the drift form by H\"older's inequality. In the two-dimensional case, the same argument with an $L^2$-drift would require an embedding of $H^{1,2}_0(U)$ into $L^\infty(U)$, which does not hold. In fact, a general $L^2$-drift need not define a bounded bilinear form on $H^{1,2}_0(U)\times H^{1,2}_0(U)$, as illustrated in \cite[Remark~2.4]{FilonovShilkin18}. The examples in \cite[Section~1.1]{Filonov13} also show that, for a general $L^2$-drift in dimension two, a solution in $H^{1,2}$ can be unbounded, and even a bounded solution need not be H\"older continuous. In particular, the first example gives a nonzero solution in $H^{1,2}_0$ of a homogeneous Dirichlet problem, so uniqueness in $H^{1,2}$ can fail.\\
Even when the drift form is bounded, the full bilinear form need not be coercive. In \cite{Trudinger73}, the weak maximum principle and Fredholm theory were used to obtain existence and uniqueness without assuming that the full bilinear form is coercive. Non-coercive linear elliptic boundary value problems with convection terms were studied in \cite{Droniou02}, where existence, uniqueness, and regularity results were established. In dimensions $d\ge3$, invertibility results and higher integrability estimates for solutions of elliptic equations with lower-order coefficients in Marcinkiewicz spaces were obtained in \cite{GGM13}.\\
For the Laplacian on bounded $C^1$ domains, existence and uniqueness of weak solutions with $H^{1,p}$-estimates were established in \cite{KimKim15} under integrability and structural assumptions on the drift. For variable diffusion coefficients, $L^p$-resolvent estimates for the Dirichlet problem were established in \cite{KangKim19} under a small-mean-oscillation assumption on the diffusion coefficients and a small Lipschitz constant for the domain, with drifts in $L^d$ when $d\ge3$ and in $L^r$ for some $r>2$ when $d=2$. Related $H^{1,p}$-estimates and solvability results for Dirichlet and Neumann problems were obtained in \cite{KangKim17}. Interior Sobolev estimates were established in \cite{Shaposhnikov08} for elliptic equations with symmetric, H\"older continuous, uniformly elliptic diffusion matrices under integrability assumptions on the lower-order coefficients. For equations with divergence-free drifts, existence, uniqueness, and regularity of weak solutions were studied in \cite{FilonovShilkin18}. In dimensions $d\ge3$, existence, uniqueness, and regularity of weak solutions for equations with the Laplacian and drifts in weak-$L^d$ were studied in \cite{KimTsai20,KimPhanTsai26} under structural assumptions on the drift. For drifts $\mathbf H\in L^{d,\infty}(U,\mathbb R^d)$ satisfying $\operatorname{div}\mathbf H\le0$ in the sense of distributions, existence and uniqueness of weak solutions in $H^{1,2+\varepsilon}_0(U)$ were established in \cite{KwonJMAA} on bounded Lipschitz domains in dimensions $d\ge2$, with bounded measurable uniformly elliptic diffusion coefficients. In \cite{L24PAMS}, the present author proved existence and uniqueness of weak solutions in $H^{1,2}_0(U)$ for equations with $L^2$-drifts satisfying $\operatorname{div}\mathbf H\le0$ in the sense of distributions, under suitable integrability assumptions on the positive zero-order coefficient and the right-hand side.\\
For a constant $c\equiv\mu\ge0$, existence and uniqueness of weak solutions in $H^{1,p}_0(U)$ were established in \cite{Kwon} for each fixed $p\in(2,\infty)$ under small mean oscillation assumptions on the diffusion coefficients and a small Lipschitz constant for the domain, together with $\operatorname{div}A\in L^2(U,\mathbb R^2)$. The proof of uniqueness is based on converting the equation into non-divergence form and applying regularity estimates and a generalized Aleksandrov maximum principle. We discuss the relation between these assumptions and our uniqueness result in Remark~\ref{comparisonkwon}.\\
In \cite{L25JMAA}, the present author studied bounded weak solutions for equations with general drifts $\mathbf H\in L^{d+\varepsilon}(U,\mathbb R^d)$ for some $\varepsilon>0$ and positive zero-order coefficients $c\in L^1(U)$. For right-hand sides $f\in L^r(U)$ with $r \in (\frac{d}{2}, \infty)$, existence and uniqueness were established by a divergence-free transformation. This transformation uses a strictly positive function $\rho$ satisfying
$$
\rho\mathbf B=A^T\nabla\rho+\rho\mathbf H,
\qquad
\operatorname{div}(\rho\mathbf B)=0\;\;  \text{ in the sense of distributions},
$$
and gives an equation with diffusion matrix $\rho A$ and drift $\rho\mathbf B$. In dimension two, under the assumption $\mathbf H\in L^{2+\varepsilon}(U,\mathbb R^2)$ for some $\varepsilon>0$, existence and uniqueness of weak solutions in $H^{1,2}_0(U)\cap L^\infty(U)$ follow from \cite[Theorem~3.6]{L25JMAA} for right-hand sides of the form $f-\operatorname{div}\mathbf F$, where $f\in L^q(U)$ and $\mathbf F\in L^q(U,\mathbb R^2)$ with $q\in(2,\infty)$. By the Hahn-Banach theorem and the duality of Lebesgue spaces, every $\ell\in H^{-1,q}(U)$ admits a representation
$$
\ell=f-\operatorname{div}\mathbf F,
\qquad
f\in L^q(U),
\qquad
\mathbf F\in L^q(U,\mathbb R^2).
$$
Consequently, the same existence and uniqueness result holds for every right-hand side $\ell\in H^{-1,q}(U)$. In dimensions $d\ge3$, this transformation was also used in \cite[Theorem~1.1]{L25ERA} to establish existence and uniqueness of bounded weak solutions in $H^{1,2}_0(U)$ for $L^d$-drifts under a divergence-type decomposition condition, with positive $c\in L^1(U)$ and right-hand sides $f\in L^q(U)$ for some $q>d/2$.\\
The uniqueness assertion in Theorem~\ref{thmmain}{\rm(i)} follows from a weak maximum principle under the assumptions $\mathbf H\in L^2(U,\mathbb R^2)$ and $c\in L^1(U)$ with $c\ge0$. More precisely, Theorem~\ref{thmwmp} shows that every $H^{1,p}_0(U)$-subsolution, with $p\in(2,\infty)$, of the homogeneous equation associated with \eqref{introproblem} satisfies $u\le0$ a.e. in $U$. The uniqueness of $H^{1,p}$-weak solutions then follows from Corollary~\ref{coruniqueness}. The diffusion matrix $A$ is only assumed to be bounded, measurable, and uniformly elliptic. No differentiability or mean oscillation condition on its entries is imposed, and no smallness or divergence condition on $\mathbf H$ is required. The result holds on every bounded open subset of $\mathbb R^2$, without any regularity assumption on the boundary.\\
Our proof uses the logarithmic construction in \cite[Theorem~3.1]{Trudinger73} and \cite[Theorem~1, p.~294]{Trudinger77}. In these arguments, a positive supremum is first shown to be attained on a set of positive measure, and a logarithmic integrability estimate gives a contradiction. In the present setting, we obtain a contradiction by examining the solution near a single maximum point. Let $w$ be the positive part of a subsolution, and suppose that its maximum is $K>0$. For $\varepsilon>0$, set
$$
v_\varepsilon:=\ln\frac{K+\varepsilon}{K+\varepsilon-w}.
$$
By choosing a suitable test function and using $c\ge0$, we obtain a bound for $\|\nabla v_\varepsilon\|_{L^2(U)}$ independent of $\varepsilon$. Passing to a subsequence as $\varepsilon$ tends to zero, we obtain $v\in H^{1,2}_0(U)$ satisfying
$$
v=\ln\frac{K}{K-w}
\qquad\text{a.e. in }U.
$$
Since $p\in(2,\infty)$, the function $w$ is H\"older continuous. Hence, if $w(x_0)=K$, then
$$
0\le K-w(x)\le C\|x-x_0\|^\alpha,
\qquad
\alpha:=1-\frac2p>0, \quad \text{ near $x_0$.}
$$
It follows that
$$
v(x)\ge\alpha\ln\frac1{\|x-x_0\|}-C_0 \quad \text{ for a.e. $x$ sufficiently close to $x_0$.}
$$
Averaging this inequality over circles centered at $x_0$ gives a logarithmic lower bound for the circular averages of $v$. By Lemma~\ref{circularineq}, these averages can grow at most at the rate $\sqrt{\ln\frac{R}{r}}$ as the radius $r$ tends to zero, since $v\in H^{1,2}_0(U)$. This contradicts the lower bound above. The assumption $p\in(2,\infty)$ is used here to ensure the H\"older continuity of $w$.\\
We next prove existence on a bounded open set with Lipschitz boundary under the additional assumption $c\in L^s(U)$ with $s>1$. The higher integrability of gradients for uniformly elliptic divergence-form equations with bounded measurable coefficients was established in \cite{Meyers63}. For the global Dirichlet problem, we use the regularity result in \cite{Groger89}. As stated in Proposition~\ref{propgroger}, there exists $\varepsilon_0>0$, depending only on the ellipticity and boundedness constants of $A$ and the geometry of $U$, such that
$$
\mathcal A_p:H^{1,p}_0(U)\longrightarrow H^{-1,p}(U),
\qquad
\mathcal A_pu:=-\operatorname{div}(A\nabla u),
$$
is a bounded isomorphism whenever $2<p<2+\varepsilon_0$. We also obtain that the drift operator $\mathcal D_{\mathbf H}$ is bounded from $H^{1,p}_0(U)$ into $H^{-1,p}(U)$. For these exponents, we further assume that
$$
s\ge\frac{2p}{p+2}.
$$
Under this condition, 
the multiplication operator $\mathcal Z_c$ is bounded from $H^{1,p}_0(U)$ into $H^{-1,p}(U)$. By approximating $\mathbf H$ by smooth functions and $c$ by bounded functions, we obtain $\varepsilon$-inequalities in which the norm of each lower-order term in $H^{-1,p}(U)$ is bounded by an arbitrarily small multiple of $\|u\|_{H^{1,p}(U)}$ and a constant multiple of $\|u\|_{L^p(U)}$. For a similar drift estimate, see also \cite[Proposition~3.1]{Kwon}. Combining these inequalities with the estimate for $\mathcal A_p$ and absorbing the small terms gives a uniform a priori estimate for
$$
\mathcal T_t=\mathcal A_p+t\mathcal D_{\mathbf H}+t\mathcal Z_c,
\qquad
0\le t\le1.
$$
This estimate still contains a term involving $\|u\|_{L^p(U)}$. To remove it, we argue by contradiction. If the uniform inverse estimate fails, we can choose $t_n\in[0,1]$ and $u_n\in H^{1,p}_0(U)$ such that
$$
\|u_n\|_{H^{1,p}(U)}=1,
\qquad
\lim_{n\to\infty}\|\mathcal T_{t_n}u_n\|_{H^{-1,p}(U)}=0.
$$
After passing to a subsequence, we obtain a limit $t\in[0,1]$ and a weak limit $u\in H^{1,p}_0(U)$. The compact embedding into $L^p(U)$, together with the a priori estimate, shows that $u$ cannot be zero. On the other hand, passing to the limit in the equation gives $\mathcal T_tu=0$, which contradicts the uniqueness result. This proves the uniform inverse estimate in Proposition~\ref{uniforminverse}. Since $\mathcal T_0=\mathcal A_p$ is invertible, the method of continuity then gives the invertibility of $\mathcal T_1$ and hence existence. These arguments yield the bounded isomorphism in Theorem~\ref{thmexistence}. For any prescribed $q \in (2, \infty)$, we choose an admissible solution exponent $p=2+\varepsilon<q$ and apply this isomorphism to $\ell\in H^{-1,q}(U)\subset H^{-1,p}(U)$. This proves Theorem~\ref{thmmain}{\rm(ii)}.\\
The assumption $c\in L^s(U)$ with $s>1$ cannot in general be replaced by $c\in L^1(U)$ in this existence result. In Proposition~\ref{endpointcounterexample}, we construct a positive coefficient $c\in L^1(B_1)$ and $f\in C_0^\infty(B_1)$ such that
$$
-\Delta u+cu=f
\qquad\text{in }B_1,
\qquad
u=0
\qquad\text{on }\partial B_1,
$$
has a bounded solution in $H^{1,2}_0(B_1)$ but has no solution in $H^{1,p}_0(B_1)$ for any $p\in(2,\infty)$. Thus mere $L^1$-integrability of the nonnegative zero-order coefficient does not suffice to guarantee an $H^{1,p}$-weak
solution for any $p>2$, even for smooth right-hand sides.
\\
This paper is organized as follows. Section~\ref{notationssec} introduces our notation and conventions. In Section~\ref{uniquenesec}, we prove the circular-mean estimate and the weak maximum principle, and then deduce uniqueness of $H^{1,p}$-weak solutions, proving Theorem~\ref{thmmain}{\rm(i)}. Section~\ref{existsec} consists of three subsections. The first establishes estimates for the lower-order terms. The second derives uniform estimates along a homotopy and proves existence by the method of continuity, completing the proof of Theorem~\ref{thmmain}{\rm(ii)}. The third presents the endpoint counterexample for $c\in L^1(U)$.

\section{Notation and conventions}\label{notationssec}
For $d\ge1$, we denote the Euclidean inner product and norm on $\mathbb R^d$ by $\langle\cdot,\cdot\rangle$ and $\|\cdot\|$, respectively. We also write $\|x\|$ for the Euclidean norm of a vector $x$. For $x_0\in\mathbb R^d$ and $r>0$, let
$$
B_r(x_0):=\{x\in\mathbb R^d:\|x-x_0\|<r\},
\qquad
B_r:=B_r(0).
$$
We denote Lebesgue measure by $dx$. For a Lebesgue measurable set $E$, we write $|E|:=dx(E)$ and denote its indicator function by $1_E$. Unless otherwise stated, integrals and almost-everywhere statements are understood with respect to Lebesgue measure.\\
For a real-valued function $u$, we write
$$
u^+:=\max\{u,0\}, \qquad u^-:=\max\{-u,0\}.
$$
Let $V\subset\mathbb R^d$ be open. We write $C_0^\infty(V)$ for the space of smooth functions with compact support in $V$, and $\mathcal D'(V)$ for the space of distributions on $V$.\\
For $1\le s\le\infty$, we denote by $L^s(V)$ the usual Lebesgue space of real-valued functions on $V$, and by $L^s(V,\mathbb R^d)$ the corresponding Lebesgue space of $\mathbb R^d$-valued functions on $V$. For $\mathbf F\in L^s(V,\mathbb R^d)$, its norm is defined as the $L^s(V)$-norm of the scalar-valued function $x\mapsto\|\mathbf F(x)\|$, namely,
$$
\|\mathbf F\|_{L^s(V,\mathbb R^d)}
:=
\bigl\|\,\|\mathbf F(\cdot)\|\,\bigr\|_{L^s(V)}.
$$
We omit the target space $\mathbb R^d$ from the notation for vector-valued norms when it is clear from the context.\\
For $1<s<\infty$, let
$$
s':=\frac{s}{s-1}
$$
be the conjugate exponent of $s$, with the endpoint conventions $1'=\infty$ and $\infty'=1$. For $1\le s<\infty$, we denote the first-order Sobolev space by
$$
H^{1,s}(V)=W^{1,s}(V)
:=
\left\{
u\in L^s(V):
\partial_j u\in L^s(V)
\text{ for }j=1,\ldots,d
\right\},
$$
where the derivatives are understood in the sense of distributions. We equip this space with the norm
$$
\|u\|_{H^{1,s}(V)}
:=\left(\|u\|_{L^s(V)}^s+\|\nabla u\|_{L^s(V)}^s\right)^{1/s}.
$$
We denote by $H_0^{1,s}(V)$ the closure of $C_0^\infty(V)$ in $H^{1,s}(V)$ and equip it with the norm inherited from $H^{1,s}(V)$. We write $u\in H_{loc}^{1,s}(V)$ if $u|_W\in H^{1,s}(W)$ for every open set $W$ compactly contained in $V$.\\
For a Banach space $X$, we denote its continuous dual by $X^*$. For $1<s<\infty$, we define
$$
H^{-1,s}(V):=\bigl(H_0^{1,s'}(V)\bigr)^*.
$$
The duality pairing between $H^{-1,s}(V)$ and $H_0^{1,s'}(V)$ is denoted by
$$
\langle\ell,\varphi\rangle_{H^{-1,s}(V)}, \qquad
\ell\in H^{-1,s}(V),
\quad
\varphi\in H_0^{1,s'}(V).
$$
We omit the subscript when the underlying spaces are clear. The dual norm is given by
$$
\|\ell\|_{H^{-1,s}(V)}:=\sup_{\substack{\varphi\in H_0^{1,s'}(V)\\
\|\varphi\|_{H^{1,s'}(V)}\le1}}
\left|
\langle\ell,\varphi\rangle_{H^{-1,s}(V)}
\right|.
$$
Each element of $H^{-1,s}(V)$ is identified with the distribution obtained by restricting it to $C_0^\infty(V)$. Differential operators are understood in the distributional sense unless otherwise stated. In particular, for $\mathbf F\in L_{loc}^1(V,\mathbb R^d)$, we define
$$
\langle\operatorname{div}\mathbf F,\varphi\rangle
:=
-\int_V\langle\mathbf F,\nabla\varphi\rangle\,dx,
\qquad
\varphi\in C_0^\infty(V).
$$
Thus, for $1<s<\infty$ and $\mathbf F\in L^s(V,\mathbb R^d)$, the same formula defines $\operatorname{div}\mathbf F\in H^{-1,s}(V)$ by continuity.\\
For $u\in H_0^{1,s}(V)$, its zero extension belongs to $H^{1,s}(\mathbb R^d)$. We denote this extension by the same symbol when no confusion can arise. In particular, if $U\subset\mathbb R^2$ is bounded and open and $p\in(2,\infty)$, the zero extension of every $u\in H_0^{1,p}(U)$ admits a representative in $C^{0,\alpha}(\mathbb R^2)$, where
$\alpha:=1-\frac2p$. This representative vanishes on $\mathbb R^2\setminus U$. Whenever pointwise values of such a function are considered, we use this representative. No regularity of $\partial U$ is required for this convention.\\
For Banach spaces $X$ and $Y$, we denote the space of bounded linear operators from $X$ to $Y$ by $\mathcal L(X,Y)$, and write $\mathcal L(X):=\mathcal L(X,X)$. The symbol $I$ denotes the identity operator on the space under consideration.
\section{Uniqueness} \label{uniquenesec}
The following elementary lemma is based on a standard polar-coordinate
argument for Sobolev functions. We give a direct proof in the present setting.
\begin{lem} \label{circularineq}
Let $B_R(x_0)$ be an open ball in $\mathbb R^2$ with center $x_0$ and radius $R$, and let $v\in H^{1,2}(B_R(x_0))$. Define $\omega:\mathbb R\to\mathbb R^2$ by
\begin{equation} \label{polarcordinome}
\omega(\theta):=(\cos\theta,\sin\theta), \qquad \theta \in \mathbb{R}.
\end{equation}
Then the following assertions hold.
\begin{itemize}
\item[(i)] For every $a,b$ with $0<a<b<R$, the function $\widetilde v:(a,b)\times(-\pi,\pi)\to\mathbb R$ defined by
$$
\widetilde v(s,\theta):=v(x_0+s\omega(\theta)), \qquad (s,\theta)\in(a,b)\times(-\pi,\pi),
$$
belongs to
$$
H^{1,2}((a,b)\times(-\pi,\pi)),
$$
and its weak derivatives are given by
\begin{align}
\partial_s\widetilde v(s,\theta)
&=
\left\langle
\nabla v(x_0+s\omega(\theta)),
\omega(\theta)
\right\rangle,
&&\text{for a.e. }(s,\theta)\in(a,b)\times(-\pi,\pi),
\label{radialderivative}
\\
\partial_{\theta}\widetilde v(s,\theta)
&=
\left\langle
\nabla v(x_0+s\omega(\theta)),
s(-\sin\theta,\cos\theta)
\right\rangle,
&&\text{for a.e. }(s,\theta)\in(a,b)\times(-\pi,\pi).
\label{thetaderivative}
\end{align}
\item[(ii)] For a.e. $r\in(0,R)$,
\begin{equation*} \label{squremeanfunc}
m(r):=\frac{1}{2\pi} \int_{-\pi}^{\pi} v(x_0+r\omega(\theta))\,d\theta
\end{equation*}
is well defined. Moreover,
$$
m \in H^{1,2}_{loc}((0,R)),
$$
and
$$
m'(s)=\frac{1}{2\pi}
\int_{-\pi}^{\pi}
\left\langle \nabla v(x_0+s\omega(\theta)), \omega(\theta) \right\rangle \,d\theta, \qquad  \text{ for a.e. $s\in(0,R)$.}
$$
In particular, $m$ has a locally absolutely continuous representative on $(0,R)$.
\item[(iii)] Let $m$ denote the locally absolutely continuous representative obtained in {\rm (ii)}. Then, for every $r, R_0$ with $0<r<R_0<R$,
\begin{equation*}\label{circularestimate}
|m(r)-m(R_0)|
\le
\frac{1}{\sqrt{2\pi}}
\|\nabla v\|_{L^2(B_{R_0}(x_0))}
\left(
\ln\frac{R_0}{r}
\right)^{\frac12}.
\end{equation*}
\end{itemize}
\end{lem}
\begin{proof}
(i) Fix $a,b$ with $0<a<b<R$ and define
$$
\Phi:(a,b)\times(-\pi,\pi)\to\mathbb{R}^2, \qquad \Phi(s,\theta):=x_0+s\omega(\theta), \qquad (s,\theta)\in(a,b)\times(-\pi,\pi).
$$
Then $\Phi$ is a $C^\infty$-diffeomorphism from $(a,b)\times(-\pi,\pi)$ onto
$$
A_{a,b}^*:=\left(B_b(x_0)\setminus\overline{B}_a(x_0) \right) \setminus \left\{ x_0-te_1:a<t<b \right\}, \qquad e_1:=(1,0).
$$
Indeed,
$$
D\Phi(s,\theta)
=
\begin{pmatrix}
\cos\theta & -s\sin\theta\\
\sin\theta & s\cos\theta
\end{pmatrix},
\qquad \det D\Phi(s,\theta)=s>0 \quad \text{ for all } (s,\theta)\in(a,b)\times(-\pi,\pi).
$$
The removed radial segment has two-dimensional Lebesgue measure zero. Since $a<s<b$, the change-of-variables formula also shows that composition with $\Phi$ preserves almost-everywhere equality. Thus the pullback of the Sobolev function $v$ is well defined up to a null set, and
$$
\widetilde v=v\circ\Phi
\qquad\text{a.e. on }(a,b)\times(-\pi,\pi).
$$
To prove that $\widetilde v$ belongs to $H^{1,2}((a,b)\times(-\pi,\pi))$, choose $v_n\in C^\infty(\overline{B}_R(x_0))$ such that
$$
\lim_{n \rightarrow \infty}v_n = v \qquad \text{ in }H^{1,2}(B_R(x_0)),
$$
and set
$$
\widetilde v_n:=v_n\circ\Phi \qquad \text{ on $(a,b) \times (-\pi, \pi)$}.
$$
The change of variables formula gives
\begin{align*}
\int_a^b\int_{-\pi}^{\pi}
|\widetilde v_n-\widetilde v_m|^2
\,d\theta\,ds
&=
\int_a^b\int_{-\pi}^{\pi}
\frac{
|v_n(x_0+s\omega(\theta))
-v_m(x_0+s\omega(\theta))|^2
}{s}
\,s\,d\theta\,ds
\\
&=
\int_{A_{a,b}^*}
\frac{|v_n(x)-v_m(x)|^2}
{\|x-x_0\|}
\,dx.
\end{align*}
Since $\|x-x_0\|\ge a$ for all $x \in A_{a,b}^*$, we obtain
\begin{align} \label{standacauch}
\int_a^b\int_{-\pi}^{\pi} |\widetilde v_n-\widetilde v_m|^2 \,d\theta\,ds \le \frac1a \int_{A_{a,b}^*} |v_n-v_m|^2\,dx= \frac1a \int_{B_b(x_0)\setminus\overline{B}_a(x_0)} |v_n-v_m|^2\,dx.
\end{align}
Moreover, by the chain rule,
\begin{equation} \label{partialsfo}
\partial_s\widetilde v_n(s,\theta)=\left\langle \nabla v_n(x_0+s\omega(\theta)), \omega(\theta) \right\rangle \quad \text{ for all  $(s, \theta) \in (a,b) \times (-\pi, \pi)$}.
\end{equation}
Since $\|\omega(\theta)\|=1$ for all $\theta \in (-\pi, \pi)$, the same change of variables gives
\begin{align*}
&\int_a^b\int_{-\pi}^{\pi}
|\partial_s\widetilde v_n
-\partial_s\widetilde v_m|^2
\,d\theta\,ds \le
\int_a^b\int_{-\pi}^{\pi} \|\nabla v_n(x_0+s\omega(\theta))-\nabla v_m(x_0+s\omega(\theta)) \|^2 \,d\theta\,ds \\
&=\int_{A_{a,b}^*}
\frac{\|\nabla v_n(x)-\nabla v_m(x)\|^2 }{\|x-x_0\| }
\,dx \le
\frac1a
\int_{B_b(x_0)\setminus\overline{B_a(x_0)}}
\|\nabla v_n-\nabla v_m\|^2\,dx.
\end{align*}
Differentiating $\widetilde v_n$ with respect to $\theta$, we obtain
\begin{align}
\partial_\theta\widetilde v_n(s,\theta)
&=s\left(-\partial_1v_n(x_0+s\omega(\theta))\sin\theta+\partial_2v_n(x_0+s\omega(\theta))\cos\theta\right) \nonumber \\
&=\left\langle\nabla v_n(x_0+s\omega(\theta)), s(-\sin\theta,\cos\theta) \right\rangle \quad \text{ for all  $(s, \theta) \in (a,b) \times (-\pi, \pi)$}. \label{partialthetafo}
\end{align}
Therefore,
\begin{align*}
&\int_a^b\int_{-\pi}^{\pi}
|\partial_\theta\widetilde v_n(s,\theta)-\partial_\theta\widetilde v_m(s,\theta)|^2 \,d\theta\,ds
\\
&=
\int_a^b\int_{-\pi}^{\pi} s^2 \big| -(\partial_1v_n-\partial_1v_m)(x_0+s\omega(\theta)) \sin\theta+(\partial_2v_n-\partial_2v_m) (x_0+s\omega(\theta)) \cos\theta \big|^2 \,d\theta\,ds \\
&\le
\int_a^b\int_{-\pi}^{\pi} s^2 \| \nabla v_n(x_0+s\omega(\theta)) -\nabla v_m(x_0+s\omega(\theta)) \|^2 \,d\theta\,ds \\
&=
\int_a^b\int_{-\pi}^{\pi} s \| \nabla v_n(x_0+s\omega(\theta))-\nabla v_m(x_0+s\omega(\theta)) \|^2 \,s\,d\theta\,ds \\
&=
\int_{A_{a,b}^*} \|x-x_0\| \|\nabla v_n(x)-\nabla v_m(x)\|^2 \,dx \\
&\le
b \int_{A_{a,b}^*} \|\nabla v_n(x)-\nabla v_m(x)\|^2 \,dx =b \int_{B_b(x_0)\setminus\overline{B}_a(x_0)} \|\nabla v_n-\nabla v_m\|^2 \,dx.
\end{align*}
Thus $(\widetilde v_n)_{n\ge1}$ is a Cauchy sequence in $H^{1,2}((a,b)\times(-\pi,\pi))$. Moreover, applying the estimate in \eqref{standacauch} to $v_n-v$, we obtain
\begin{align*}
\int_a^b\int_{-\pi}^{\pi} |\widetilde v_n-\widetilde v|^2 \,d\theta\,ds
&\le \frac1a \int_{B_b(x_0)\setminus\overline{B}_a(x_0)} |v_n-v|^2\,dx.
\end{align*}
Since
$$
\lim_{n\to\infty}\frac1a\int_{B_b(x_0)\setminus\overline{B}_a(x_0)}|v_n-v|^2\,dx=0,
$$
the left-hand side also has limit zero. Hence,
$$
\lim_{n \rightarrow \infty} \widetilde v_n=\widetilde v=v\circ\Phi
\qquad
\text{in }L^2((a,b)\times(-\pi,\pi)).
$$
Consequently,
$$
\widetilde v
\in
H^{1,2}((a,b)\times(-\pi,\pi))
$$
and
$$
\lim_{n \rightarrow \infty} \widetilde v_n=\widetilde v=v\circ\Phi
\qquad
\text{in } H^{1,2}((a,b)\times(-\pi,\pi)).
$$
To identify the derivatives, the same change of variables gives
$$
\int_a^b\int_{-\pi}^{\pi}
\|\nabla v_n(x_0+s\omega(\theta))
-\nabla v(x_0+s\omega(\theta))\|^2\,d\theta\,ds
\le\frac1a\|\nabla v_n-\nabla v\|_{L^2(B_R(x_0))}^2.
$$
The right-hand side has limit zero because of the $H^{1,2}(B_R(x_0))$ convergence of $v_n$. Since $\|\omega(\theta)\|=1$ and $\|s(-\sin\theta,\cos\theta)\|\le b$, the right-hand sides of \eqref{partialsfo} and \eqref{partialthetafo} therefore converge in $L^2((a,b)\times(-\pi,\pi))$ to the corresponding expressions with $v$ in place of $v_n$. Their left-hand sides converge to the weak derivatives of $\widetilde v$ by the $H^{1,2}$ convergence already proved. This gives \eqref{radialderivative} and \eqref{thetaderivative}, and completes the proof of {\rm (i)}.\\ \\
(ii) Fix $a,b$ with $0<a<b<R$. By (i), $\widetilde v\in L^2((a,b)\times(-\pi,\pi))$. Hence, by Fubini's theorem,
$$
\int_{-\pi}^{\pi}
|\widetilde v(s,\theta)|^2\,d\theta<\infty \quad \text{ for a.e. $s\in(a,b)$. }
$$
Thus, for a.e. $s \in (a,b)$, the Cauchy-Schwarz inequality gives
$$
\int_{-\pi}^{\pi}
|\widetilde v(s,\theta)|\,d\theta
\le
\sqrt{2\pi}
\left(
\int_{-\pi}^{\pi}
|\widetilde v(s,\theta)|^2\,d\theta
\right)^{1/2}
<\infty.
$$
Therefore, for a.e. $s\in(a,b)$,
$$
m(s)=\frac{1}{2\pi} \int_{-\pi}^{\pi} \widetilde v(s,\theta)\,d\theta
$$
is well defined. By Jensen's inequality, for a.e. $s\in(a,b)$,
\begin{align*}
|m(s)|^2 =\left| \frac{1}{2\pi} \int_{-\pi}^{\pi} \widetilde v(s,\theta)\,d\theta \right|^2 \le \frac{1}{2\pi} \int_{-\pi}^{\pi} |\widetilde v(s,\theta)|^2\,d\theta.
\end{align*}
Integrating over $(a,b)$ yields
$$
\int_a^b|m(s)|^2\,ds
\le
\frac{1}{2\pi}
\int_a^b\int_{-\pi}^{\pi}
|\widetilde v(s,\theta)|^2\,d\theta\,ds
<\infty.
$$
Thus, $m\in L^2(a,b)$. Let $\psi\in C_0^\infty((a,b))$. By Fubini's theorem and {\rm (i)},
\begin{align*}
\int_a^b m(s)\psi'(s)\,ds
&=
\frac{1}{2\pi}
\int_a^b\int_{-\pi}^{\pi}
\widetilde v(s,\theta)\psi'(s)
\,d\theta\,ds
\\
&=
-\frac{1}{2\pi}
\int_a^b\int_{-\pi}^{\pi}
\partial_s\widetilde v(s,\theta)\psi(s)
\,d\theta\,ds
=
-\int_a^b
\left(
\frac{1}{2\pi}
\int_{-\pi}^{\pi}
\partial_s\widetilde v(s,\theta)\,d\theta
\right)
\psi(s)\,ds.
\end{align*}
For completeness, the integration by parts in $s$ can first be applied with the test function $\psi(s)\chi_j(\theta)$, where $\chi_j\in C_0^\infty((-\pi,\pi))$, $0\le\chi_j\le1$, and $\lim_{j\to\infty}\chi_j(\theta)=1$ for every $\theta\in(-\pi,\pi)$. There is no derivative of $\chi_j$ in this identity. Since both $\widetilde v$ and $\partial_s\widetilde v$ belong to $L^2((a,b)\times(-\pi,\pi))$, dominated convergence gives the displayed identity without the cutoff. Hence the weak derivative of $m$ is
$$
m'(s)
=
\frac{1}{2\pi}
\int_{-\pi}^{\pi}
\partial_s\widetilde v(s,\theta)\,d\theta \quad \text{ for a.e. $s\in(a,b)$.}
$$
By \eqref{radialderivative},
\begin{equation} \label{mderivaform}
m'(s)=\frac{1}{2\pi} \int_{-\pi}^{\pi} \left\langle \nabla v(x_0+s\omega(\theta)), \omega(\theta) \right\rangle \,d\theta \quad \text{ for a.e. $s\in(a,b)$.}
\end{equation}
Moreover, by Jensen's inequality and the Cauchy-Schwarz inequality,
$$
|m'(s)|^2 \le \frac{1}{2\pi} \int_{-\pi}^{\pi} |\partial_s\widetilde v(s,\theta)|^2\,d\theta \quad \text{ for a.e. $s\in(a,b)$.}
$$
Therefore, $m'\in L^2\big((a,b)\big)$, and hence $m\in H^{1,2}\big((a,b)\big)$. Since $a,b$ were arbitrary with $0<a<b<R$, we conclude that
$$
m\in H^{1,2}_{loc}((0,R)).
$$
Since every one-dimensional $H^{1,2}\big((a,b)\big)$ function admits an absolutely continuous representative, $m$ admits a locally absolutely continuous representative on $(0,R)$. We denote this representative again by $m$. By \eqref{mderivaform}, it satisfies
$$
m'(s)=\frac{1}{2\pi} \int_{-\pi}^{\pi} \left\langle \nabla v(x_0+s\omega(\theta)), \omega(\theta) \right\rangle \,d\theta \quad \text{ for a.e. $s\in(0,R)$.}
$$
This proves {\rm (ii)}.\\ \\
(iii) Let $m$ denote the locally absolutely continuous representative obtained in {\rm (ii)}. Choose $r$ and $R_0$ with $0<r<R_0<R$. The fundamental theorem of calculus implies
$$
m(r)-m(R_0)=-\int_r^{R_0}m'(s)\,ds.
$$
By the formula for $m'$ in {\rm (ii)}, the Cauchy-Schwarz inequality, and $\|\omega(\theta)\|=1$,
\begin{align*}
|m'(s)| &\le \frac{1}{2\pi} \int_{-\pi}^{\pi} \|\nabla v(x_0+s\omega(\theta))\| \,d\theta \le \frac{1}{\sqrt{2\pi}} \left( \int_{-\pi}^{\pi} \|\nabla v(x_0+s\omega(\theta))\|^2 \,d\theta \right)^{1/2}.
\end{align*}
Applying the Cauchy-Schwarz inequality in the radial variable then gives
\begin{align*}
|m(r)-m(R_0)|
&\le
\int_r^{R_0}|m'(s)|\,ds
\\
&\le
\frac{1}{\sqrt{2\pi}}
\int_r^{R_0}
\left( \int_{-\pi}^{\pi} \|\nabla v(x_0+s\omega(\theta))\|^2 \,d\theta \right)^{1/2} ds \\
&=
\frac{1}{\sqrt{2\pi}}
\int_r^{R_0}
\left(
\int_{-\pi}^{\pi}
\|\nabla v(x_0+s\omega(\theta))\|^2
s\,d\theta
\right)^{1/2}
\frac{1}{\sqrt{s}}\,ds
\\
&\le
\frac{1}{\sqrt{2\pi}}
\left(
\int_r^{R_0}
\int_{-\pi}^{\pi}
\|\nabla v(x_0+s\omega(\theta))\|^2
s\,d\theta\,ds
\right)^{1/2}
\left(
\int_r^{R_0}\frac{1}{s}\,ds
\right)^{1/2}.
\end{align*}
By the change of variables formula,
\begin{align*}
&\int_r^{R_0}
\int_{-\pi}^{\pi}
\|\nabla v(x_0+s\omega(\theta))\|^2 s\,d\theta\,ds=\int_{B_{R_0}(x_0)\setminus\overline{B}_r(x_0)}
\|\nabla v(x)\|^2\,dx \le \|\nabla v\|_{L^2(B_{R_0}(x_0))}^2.
\end{align*}
Therefore,
\begin{align*}
|m(r)-m(R_0)|
&\le \frac{1}{\sqrt{2\pi}} \|\nabla v\|_{L^2(B_{R_0}(x_0))} \left( \int_r^{R_0}\frac{ds}{s}
\right)^{1/2}=\frac{1}{\sqrt{2\pi}} \|\nabla v\|_{L^2(B_{R_0}(x_0))} \left(\ln\frac{R_0}{r} \right)^{\frac12},
\end{align*}
as desired.
\end{proof}
\\
We now establish the weak maximum principle under the assumption that $\mathbf H\in L^2(U,\mathbb R^2)$ and $c\in L^1(U)$ with $c\ge0$ in $U$, without any smallness, sign, or divergence condition on the drift. The proof uses the logarithmic construction in \cite[Theorem~3.1]{Trudinger73} and \cite[Theorem~1]{Trudinger77}, followed by a contradiction at a single maximum point. If a subsolution in $H^{1,p}_0(U)$, with $p\in(2,\infty)$, attains a positive maximum, its H\"older continuity forces the circular mean of a logarithmic transform to grow at least as fast as $\ln\frac{1}{r}$. Lemma~\ref{circularineq}, however, bounds the same mean by a constant multiple of $\sqrt{\ln\frac{1}{r}}$. This is the final step of the proof below.
\begin{theo}[Weak maximum principle]\label{thmwmp}
Assume that $U$ is a bounded open subset of $\mathbb R^2$,
$\mathbf H\in L^2(U,\mathbb R^2)$, and $c\in L^1(U)$ with
$c\ge0$ a.e. in $U$. Assume also that the matrix
$A=(a_{ij})_{1\le i,j\le2}$ has bounded measurable entries on
$\mathbb R^2$ and satisfies \eqref{elliptbdda}.
Let $p\in(2,\infty)$ and suppose that
$u\in H^{1,p}_0(U)\cap C(\overline U)$ satisfies
\begin{equation}\label{subsolutionequat}
\int_U \langle A\nabla u,\nabla\varphi\rangle\,dx +
\int_U \langle\mathbf{H},\nabla u\rangle\varphi\,dx + \int_U
cu\varphi\,dx \le 0 \quad \text{ for every $\varphi\in C_0^\infty(U)$ with $\varphi \geq 0$}.
\end{equation}
Then
$$
u(x)\le0 \qquad\text{ for all $x \in \overline{U}$}.
$$
\end{theo}
\begin{proof}
We identify $u$ with its zero extension to $\mathbb R^2$. Since $u\in H^{1,p}_0(U)$, its zero extension belongs to $H^{1,p}(\mathbb R^2)$. As $p\in(2,\infty)$, Morrey's inequality yields a continuous representative of this extension. We use this representative throughout. In particular,
$$
u=0 \qquad \text{on }\mathbb R^2\setminus U.
$$
For fixed $u\in H^{1,p}_0(U)$, each of the three terms in \eqref{subsolutionequat} is continuous with respect to the $H^{1,p}_0(U)$-norm of the test function. Indeed, since $U$ is bounded,
\begin{align*}
\left|
\int_U
\langle A\nabla u,\nabla\varphi\rangle\,dx
\right|
\le
2M
\int_U
\|\nabla u\|\,\|\nabla\varphi\|\,dx \le
2M|U|^{1-\frac{2}{p}}
\|\nabla u\|_{L^p(U)}
\|\nabla\varphi\|_{L^p(U)}.
\end{align*}
Moreover, setting $q:=\frac{2p}{p-2}$, we have $\frac12+\frac1p+\frac1q=1$, and hence
\begin{align*}
\left|
\int_U
\langle\mathbf H,\nabla u\rangle\varphi\,dx
\right| \le \|\mathbf H\|_{L^2(U)} \|\nabla u\|_{L^p(U)} \|\varphi\|_{L^q(U)} \le C_1\|\mathbf H\|_{L^2(U)} \|u\|_{H^{1,p}_0(U)} \|\varphi\|_{H^{1,p}_0(U)},
\end{align*}
where $C_1>0$ depends only on $p$ and $|U|$. Finally,
\begin{align*}
\left|
\int_U cu\varphi\,dx
\right| \le \|c\|_{L^1(U)} \|u\|_{L^\infty(U)}
\|\varphi\|_{L^\infty(U)} \le C_2\|c\|_{L^1(U)}
\|u\|_{H^{1,p}_0(U)}
\|\varphi\|_{H^{1,p}_0(U)}.
\end{align*}
We next justify the use of nonnegative test functions in $H^{1,p}_0(U)$. Let $\varphi\in H^{1,p}_0(U)$ with $\varphi\ge0$, and choose $\psi_n\in C_0^\infty(U)$ such that
$$
\lim_{n\to\infty}\psi_n=\varphi
\qquad\text{in }H^{1,p}(U).
$$
Continuity of the positive-part map in $H^{1,p}(U)$ gives
$$
\lim_{n\to\infty}\psi_n^+=\varphi
\qquad\text{in }H^{1,p}(U).
$$
Each $\psi_n^+$ is nonnegative and has compact support in $U$. Mollifying it with a nonnegative smooth kernel of sufficiently small support, we may choose $\varphi_n\in C_0^\infty(U)$ such that
$$
\varphi_n\ge0,
\qquad
\|\varphi_n-\psi_n^+\|_{H^{1,p}(U)}<\frac1n.
$$
The radius of mollification is chosen smaller than the distance from $\operatorname{supp}\psi_n$ to $\mathbb R^2\setminus U$, so $\varphi_n$ remains compactly supported in $U$. It follows that $\lim_{n\to\infty}\varphi_n=\varphi$ in $H^{1,p}_0(U)$. Applying \eqref{subsolutionequat} to $\varphi_n$ and using the three continuity estimates above proves \eqref{subsolutionequat} for $\varphi$. Thus every nonnegative member of $H^{1,p}_0(U)$ is an admissible test function.\\
Now set
$$
w:=u^+ \qquad \text{ in $\mathbb{R}^2$}.
$$
The positive-part map preserves $H^{1,p}_0(U)$, so the restriction of $w$ to $U$ belongs to this space. Its zero extension satisfies
$$
w\in H^{1,p}(\mathbb{R}^2)\cap C(\mathbb{R}^2),
\qquad
w=0 \quad \text{on }\mathbb{R}^2\setminus U,
$$
and
$$
\nabla w= 1_{\{u>0\}}\nabla u \qquad \text{a.e. in }\mathbb{R}^2.
$$
By Morrey's inequality, there exists a constant $C_3>0$ such that
\begin{equation}\label{holderwes}
|w(x)-w(y)| \le C_3\|x-y\|^\alpha \qquad\text{for all }x,y\in\mathbb R^2,
\end{equation}
where $\alpha=1-\frac2p$. Since $U$ is bounded and $w=0$ on $\mathbb R^2\setminus U$, the function $w$ has compact support. Hence $w$ attains its maximum on $\mathbb R^2$. Suppose, for a contradiction, that
$$
K:=\max_{\mathbb{R}^2}w>0.
$$
Since $w=0$ on $\mathbb{R}^2\setminus U$, there exists $x_0\in U$ such that $w(x_0)=K$. For $\varepsilon>0$, define
$$
\Phi_\varepsilon(t):=\frac{1}{K+\varepsilon-t}-\frac{1}{K+\varepsilon}, \qquad 0\le t\le K.
$$
Then
$$
\Phi_\varepsilon(0)=0, \qquad \Phi_\varepsilon'(t) = \frac{1}{(K+\varepsilon-t)^2}>0.
$$
For each fixed $\varepsilon>0$, the function $\Phi_\varepsilon$ is smooth on $[0,K]$, vanishes at $0$, and satisfies $0<\Phi_\varepsilon'(t)\le\varepsilon^{-2}$ for all $t \in [0,K]$. Since $0\le w\le K$, the Sobolev chain rule and the zero boundary value of $w$ give
$$
\varphi_\varepsilon:=\Phi_\varepsilon(w) \in H^{1,p}_0(U), \qquad \varphi_\varepsilon\ge0,
$$
and
$$
\nabla\varphi_\varepsilon=\Phi_\varepsilon'(w)\nabla w.
$$
Notice that
$$
\varphi_\varepsilon=0 \qquad \text{ on }\{u\le0\},
$$
whereas
$$
u=w, \qquad \nabla u=\nabla w \qquad \text{ on $\{u>0\}$}.
$$
Therefore, using $\varphi_\varepsilon$ as a test function in \eqref{subsolutionequat}, we obtain
\begin{align}
&
\int_U \Phi_\varepsilon'(w) \langle A\nabla w,\nabla w\rangle\,dx + \int_U \langle\mathbf H,\nabla w\rangle \Phi_\varepsilon(w)\,dx  + \int_U cw\Phi_\varepsilon(w)\,dx \le0. \label{logtestestim}
\end{align}
Since
$$
c\ge0, \qquad w\ge0, \qquad \Phi_\varepsilon(w)\ge0, \quad \text{ in $U$},
$$
$\int_U cw\Phi_\varepsilon(w)\,dx \ge0$.  Applying the ellipticity assumption to \eqref{logtestestim}, we obtain
\begin{equation}\label{logprestima}
\lambda \int_U \frac{\|\nabla w\|^2} {(K+\varepsilon-w)^2}\,dx \le \left|
\int_U \langle\mathbf H,\nabla w\rangle \Phi_\varepsilon(w)\,dx \right|.
\end{equation}
Define
$$
v_\varepsilon:=\ln\frac{K+\varepsilon}{K+\varepsilon-w} \quad \text{ on $\mathbb{R}^2$}.
$$
Since the function
$$
t\mapsto \ln\frac{K+\varepsilon}{K+\varepsilon-t}
$$
is smooth on $[0,K]$, vanishes at $t=0$, and has derivative bounded by $\varepsilon^{-1}$ on this interval, the Sobolev chain rule gives $v_\varepsilon\in H^{1,p}_0(U)\cap C(\overline U)$. Furthermore,
$$
\nabla v_\varepsilon = \frac{\nabla w}{K+\varepsilon-w}, \qquad \Phi_\varepsilon(w) = \frac{w}{(K+\varepsilon)(K+\varepsilon-w)} \quad \text{ in $U$}.
$$
Consequently,
$$
\langle\mathbf H,\nabla w\rangle \Phi_\varepsilon(w) = \frac{w}{K+\varepsilon} \langle\mathbf H,\nabla v_\varepsilon\rangle \qquad \text{in }U.
$$
Since
$$
0\le\frac{w}{K+\varepsilon}\le1,
$$
it follows from \eqref{logprestima} and the Cauchy-Schwarz inequality that
$$
\lambda \|\nabla v_\varepsilon\|_{L^2(U)}^2 \le \|\mathbf H\|_{L^2(U)} \|\nabla v_\varepsilon\|_{L^2(U)}.
$$
Hence
\begin{equation}\label{uniformlogest}
\|\nabla v_\varepsilon\|_{L^2(U)} \le \frac{1}{\lambda} \|\mathbf H\|_{L^2(U)}.
\end{equation}
Since $U$ is bounded and $p\in(2,\infty)$, we have $H^{1,p}_0(U)\subset H^{1,2}_0(U)$. Consequently, Poincar\'e's inequality and \eqref{uniformlogest} yield
$$
\sup_{\varepsilon>0} \|v_\varepsilon\|_{H^{1,2}_0(U)}<\infty.
$$
By weak compactness, there exists a decreasing sequence $(\varepsilon_k)_{k \geq 1}$ with $\lim_{k \rightarrow \infty} \varepsilon_{k}=0$ such that
$$
\lim_{k \rightarrow \infty}v_{\varepsilon_k} =v \qquad \text{weakly in }H^{1,2}_0(U).
$$
The embedding $H^{1,2}_0(U)\hookrightarrow L^2(U)$ is compact even though no regularity of $\partial U$ is assumed. Indeed, the zero extensions belong to $H^{1,2}_0(B)$ for a fixed ball $B$ containing $\overline U$, and the Rellich-Kondrachov theorem on $B$ applies to these extensions. Thus, passing to a further subsequence if necessary, we may also assume that
$$
\lim_{k \rightarrow \infty}v_{\varepsilon_k}= v \qquad \text{ in }L^2(U).
$$
Passing to another subsequence, we also obtain
$$
\lim_{k \rightarrow \infty} v_{\varepsilon_k}(x)= v(x) \qquad \text{for a.e. }x\in U.
$$
Let $E\subset U$ be a set of full measure on which $\lim_{k\to\infty}v_{\varepsilon_k}(x)=v(x)$ and $v(x)$ is finite. If $x\in E$ and $w(x)=K$, then
$$
\lim_{k\to\infty}v_{\varepsilon_k}(x)
=
\lim_{k\to\infty}\ln\frac{K+\varepsilon_k}{\varepsilon_k}
=\infty,
$$
which is impossible. Hence
$$
|\{x\in U:w(x)=K\}|=0.
$$
Therefore, for a.e. $x\in U$, we obtain $w(x)<K$ and
$$
\lim_{k\to\infty}v_{\varepsilon_k}(x)=\ln\frac{K}{K-w(x)}.
$$
Hence
\begin{equation}\label{vlimitest}
v(x)=\ln\frac{K}{K-w(x)} \qquad \text{ for a.e. } x\in U.
\end{equation}
By \eqref{holderwes},
$$
0\le K-w(x) = w(x_0)-w(x) \le C_3\|x-x_0\|^\alpha \quad \text{ for every $x\in\mathbb{R}^2$. }
$$
Choose $R>0$ such that $\overline{B}_R(x_0)\subset U$. Using \eqref{vlimitest}, we obtain
\begin{equation}\label{vlowereq}
v(x)
\ge
\alpha\ln\frac{1}{\|x-x_0\|}-C_0
\qquad
\text{for a.e. $x\in B_R(x_0)$,}
\end{equation}
where $C_0:=1+\left|\ln\frac K{C_3}\right|>0$ is a constant independent of $x$. On the other hand, applying Lemma~\ref{circularineq}{\rm (ii)} to $v|_{B_R(x_0)}$, for $\omega$ defined in \eqref{polarcordinome}, we see that, for a.e. $r\in(0,R)$,
$$
m(r):=\frac{1}{2\pi} \int_{-\pi}^{\pi} v(x_0+r\omega(\theta))\,d\theta
$$
is well defined, and $m\in H^{1,2}_{loc}((0,R))$. Let $m$ denote the locally absolutely continuous representative given by Lemma~\ref{circularineq}{\rm (ii)}. We next derive a lower bound for $m$. Let
$$
N:=\left\{x\in B_R(x_0)\setminus\{x_0\}:
 v(x)<\alpha\ln\frac1{\|x-x_0\|}-C_0\right\}.
$$
By \eqref{vlowereq}, $|N|=0$. Hence, by the polar-coordinate formula and Fubini's theorem,
$$
0=\int_{B_R(x_0)} 1_N(x)\,dx=\int_0^R r \int_{-\pi}^{\pi} 1_N(x_0+r\omega(\theta))\,d\theta\,dr.
$$
Thus, for a.e. $r\in(0,R)$,
$$
1_N(x_0+r\omega(\theta))=0 \qquad \text{for a.e. $\theta\in(-\pi,\pi)$.}
$$
Consequently, for a.e. $r\in(0,R)$,
$$
v(x_0+r\omega(\theta)) \ge \alpha\ln\frac{1}{\|r\omega(\theta)\|}-C_0=\alpha\ln\frac1r-C_0 \qquad \text{for a.e. $\theta\in(-\pi,\pi)$.}
$$
Averaging this inequality with respect to $\theta$ over $(-\pi, \pi)$ yields
\begin{equation}\label{meanlowerest}
m(r) \ge \alpha\ln\frac1r-C_0 \qquad \text{for a.e. $r\in(0,R)$.}
\end{equation}
Fix $R_0\in(0,R)$. By Lemma~\ref{circularineq}{\rm (iii)},
\begin{equation}\label{meanupperineq}
|m(r)| \le |m(R_0)| +\frac{1}{\sqrt{2\pi}} \|\nabla v\|_{L^2(B_{R_0}(x_0))} \left( \ln\frac{R_0}{r} \right)^{\frac12}, \qquad \text{ for every $r \in (0, R_0)$.}
\end{equation}
Choose a decreasing sequence $(r_n)_{n\ge1}$ from the full-measure subset of $(0,R_0)$ on which \eqref{meanlowerest} holds such that $\lim_{n\to\infty}r_n=0$, and set
$$
t_n:= \ln \frac{R_0}{r_n}, \quad n \geq 1.
$$
Since
$$
\ln\frac1{r_n} = t_n-\ln R_0 \quad \text{ for all $n \geq 1$},
$$
combining \eqref{meanlowerest} and \eqref{meanupperineq} gives
$$
\alpha t_n \le C_1+C_2 t_n^{\frac12} \quad \text{ for all $n \geq 1$},
$$
where $C_1=C_0+|m(R_0)|+\alpha |\ln R_0|$ and $C_2=\frac{1}{\sqrt{2\pi}} \|\nabla v\|_{L^2(B_{R_0}(x_0))}$. Since $\lim_{n\to\infty}t_n=\infty$, division by $t_n$ and passage to the limit yield
$$
\alpha \le \lim_{n\to\infty}\left(\frac{C_1}{t_n}+\frac{C_2}{t_n^{\frac12}}\right)=0,
$$
which contradicts
$$
\alpha=1-\frac2p>0.
$$
Therefore the assumption $K=\max_{\mathbb{R}^2}w>0$ is false. Since $w=u^+\ge0$, we have
$$
w(x)=0 \qquad \text{for every }x\in \mathbb{R}^2.
$$
Thus
$$
u(x)\le0
\qquad
\text{for every }x\in \overline{U}.
$$
\end{proof}
\begin{cor}[Uniqueness]\label{coruniqueness}
Assume that $U$, $A$, $\mathbf H$, and $c$ satisfy the assumptions
of Theorem~\ref{thmwmp}, and let $p\in(2,\infty)$. Suppose that $u\in H^{1,p}_0(U)$ satisfies
\begin{equation}\label{homogeneouseq}
\int_U \langle A\nabla u,\nabla\varphi\rangle\,dx
+
\int_U \langle\mathbf H,\nabla u\rangle\varphi\,dx
+
\int_U cu\varphi\,dx
=
0
\end{equation}
for every $\varphi\in C_0^\infty(U)$. Then
$$
u=0
\qquad
\text{a.e. in }U.
$$
\end{cor}
\begin{proof}
Let $\widetilde u$ be the continuous representative of $u$ on $\overline U$. Since \eqref{homogeneouseq} holds with equality, $\widetilde u$ satisfies the hypotheses of Theorem~\ref{thmwmp}. Hence,
$$
\widetilde u(x)\le0
\qquad
\text{for every }x\in\overline U.
$$
By linearity, $-\widetilde u$ also satisfies the same homogeneous equation, and thus another application of Theorem~\ref{thmwmp} gives
$$
-\widetilde u(x)\le0
\qquad
\text{for every }x\in\overline U.
$$
Therefore,
$$
\widetilde u(x)=0
\qquad
\text{for every }x\in\overline U,
$$
and hence $u=0$ a.e. in $U$.
\end{proof}
\\ \\
\begin{proof}{\bf of Theorem~\ref{thmmain}~{\rm(i)}.}
Fix $p\in(2,\infty)$ and $\ell\in\mathcal D'(U)$. Let $u_1,u_2\in H^{1,p}_0(U)$ be weak solutions of \eqref{introproblem} with the same right-hand side $\ell$, and set $w:=u_1-u_2$. Then $w\in H^{1,p}_0(U)$. Subtracting the identities in \eqref{introweakformulation} for $u_1$ and $u_2$, we obtain
$$
\int_U\langle A\nabla w,\nabla\varphi\rangle\,dx
+\int_U\langle\mathbf H,\nabla w\rangle\varphi\,dx
+\int_Ucw\varphi\,dx=0
\qquad
\text{for every }\varphi\in C_0^\infty(U).
$$
The hypotheses of Corollary~\ref{coruniqueness} are satisfied, so $w=0$ a.e. in $U$. Hence $u_1=u_2$ a.e. in $U$. In particular, taking $\ell=0$ and $u_2=0$ proves part {\rm(i)}.
\end{proof}
\begin{rem}\label{comparisonkwon}
In \cite[Theorem~4.1]{Kwon}, uniqueness was proved for a constant $c\equiv\mu\ge0$ under a small-mean-oscillation assumption on $A$ and the additional condition $\operatorname{div}A\in L^2(U,\mathbb R^2)$. The proof is based on converting the equation from divergence form to non-divergence form. Krylov's solvability theory is then used to obtain local $H^{2,q}$ regularity for a suitable $q\in(1,2)$, and uniqueness follows from Krylov's generalized Aleksandrov maximum principle.\\
Taking $c\equiv\mu\ge0$ in Corollary~\ref{coruniqueness}, we obtain the same uniqueness result without the condition $\operatorname{div}A\in L^2(U,\mathbb R^2)$ or the small-mean-oscillation assumption on $A$. The possibility of removing the condition $\operatorname{div}A\in L^2(U,\mathbb R^2)$ from the existence and uniqueness result was raised in \cite[Remark~{\rm(iii)} following Theorem~2.5]{Kwon}. Indeed, the assumption $\operatorname{div}A\in L^2(U,\mathbb R^2)$ in \cite[Theorem~2.5]{Kwon} can also be removed, with all other assumptions retained, by combining Corollary~\ref{coruniqueness} with the solvability result for sufficiently large $\mu$ in \cite[Theorem~3.3]{Kwon} and applying the Fredholm argument in \cite[Section~5]{Kwon}.
\end{rem}

\section{Existence}\label{existsec}
\subsection{Estimates for lower-order terms}
For the estimates below, we assume that $c\in L^s(U)$ for some $1<s<\infty$, while retaining the assumption $\mathbf H\in L^2(U,\mathbb R^2)$. We establish boundedness and $\varepsilon$-inequalities for the multiplication operator $\mathcal Z_c$ and the drift operator $\mathcal D_{\mathbf H}$ under the conditions on the exponents specified below. The boundedness estimates follow from H\"older's inequality and the two-dimensional Sobolev embeddings. The $\varepsilon$-inequalities are of Gerhardt type (see \cite{Gerhardt79}). They are obtained by approximating $c$ by bounded functions and $\mathbf H$ by smooth functions, together with integration by parts for the drift term. Each operator is thus decomposed into a part with arbitrarily small norm from $H^{1,p}_0(U)$ to $H^{-1,p}(U)$ and a remainder controlled by the $L^p(U)$-norm of $u$.
\begin{prop}[boundedness and $\varepsilon$-inequality for $\mathcal{Z}_c$]\label{propzeroorder}
Let $p\in(2,\infty)$ and $s\in(1,\infty)$. Assume that
\begin{equation}\label{spcondition}
s\ge\frac{2p}{p+2}
\end{equation}
and $c\in L^s(U)$. Then multiplication by $c$ defines a bounded operator
\begin{equation} \label{zeroordermult}
\mathcal{Z}_c:H^{1,p}_0(U)\longrightarrow H^{-1,p}(U),
\qquad
\langle\mathcal{Z}_cu,\varphi\rangle:=\int_U cu\varphi\,dx,
\end{equation}
for $u\in H^{1,p}_0(U)$ and $\varphi\in H^{1,p'}_0(U)$, where $p':=\frac{p}{p-1}$. More precisely, there exists a constant $C_Z>0$, depending only on $p,s,|U|$, such that
\begin{equation}\label{mcbounde}
\|\mathcal{Z}_cu\|_{H^{-1,p}(U)}
\le C_Z\|c\|_{L^s(U)}\|u\|_{H^{1,p}(U)}
\qquad\text{for all }u\in H^{1,p}_0(U).
\end{equation}
Moreover, for every constant $\eta>0$, choose $c_\eta\in L^\infty(U)$ such that
$$
C_Z\|c-c_\eta\|_{L^s(U)}\le\eta,
$$
where $C_Z$ is the constant in \eqref{mcbounde}. Then
\begin{equation}\label{cepsilon}
\|\mathcal{Z}_cu\|_{H^{-1,p}(U)}
\le \eta\|u\|_{H^{1,p}(U)}
+\|c_\eta\|_{L^\infty(U)}\|u\|_{L^p(U)}
\qquad\text{for all }u\in H^{1,p}_0(U).
\end{equation}
\end{prop}
\begin{proof}
Since $p\in(2,\infty)$, we have $1<p'<2$, and the two-dimensional Sobolev embedding gives
\begin{equation} \label{bddetima}
H^{1,p'}_0(U)\hookrightarrow L^{\frac{2p}{p-2}}(U).
\end{equation}
Condition \eqref{spcondition} is equivalent to
$$
s'\le\frac{2p}{p-2},
\qquad
s':=\frac{s}{s-1}.
$$
Since $U$ is bounded, it follows from \eqref{bddetima} that
$$
H^{1,p'}_0(U)\hookrightarrow L^{s'}(U).
$$
On the other hand, since $p\in(2,\infty)$, Morrey's inequality gives
$$
H^{1,p}_0(U)\hookrightarrow L^\infty(U).
$$
Therefore, for every $u\in H^{1,p}_0(U)$ and $\varphi\in H^{1,p'}_0(U)$,
\begin{align*}
\left|\int_U cu\varphi\,dx\right|
&\le \|c\|_{L^s(U)}\|u\|_{L^\infty(U)}\|\varphi\|_{L^{s'}(U)} \le C_Z\|c\|_{L^s(U)}\|u\|_{H^{1,p}(U)}
\|\varphi\|_{H^{1,p'}(U)}.
\end{align*}
This proves \eqref{mcbounde}.\\
Let $\eta>0$ be a fixed constant. Since $s<\infty$, bounded functions are dense in $L^s(U)$, so $c_\eta$ can be chosen as in the statement. Applying \eqref{mcbounde} to $c-c_\eta$, we obtain
$$
\|(\mathcal{Z}_c-\mathcal{Z}_{c_\eta})u\|_{H^{-1,p}(U)}
\le C_Z\|c-c_\eta\|_{L^s(U)}\|u\|_{H^{1,p}(U)}
\le \eta\|u\|_{H^{1,p}(U)}
\qquad\text{for all }u\in H^{1,p}_0(U).
$$
Moreover, for every $\varphi\in H^{1,p'}_0(U)$,
$$
\left|\int_U c_\eta u\varphi\,dx\right|
\le \|c_\eta\|_{L^\infty(U)}\|u\|_{L^p(U)}
\|\varphi\|_{L^{p'}(U)}
\le \|c_\eta\|_{L^\infty(U)}\|u\|_{L^p(U)}
\|\varphi\|_{H^{1,p'}(U)}.
$$
Hence
\begin{align*}
\|\mathcal{Z}_cu\|_{H^{-1,p}(U)}
&\le \|(\mathcal{Z}_c-\mathcal{Z}_{c_\eta})u\|_{H^{-1,p}(U)}
+\|\mathcal{Z}_{c_\eta}u\|_{H^{-1,p}(U)}\\
&\le \eta\|u\|_{H^{1,p}(U)}
+\|c_\eta\|_{L^\infty(U)}\|u\|_{L^p(U)}
\qquad\text{for all }u\in H^{1,p}_0(U),
\end{align*}
as desired.
\end{proof}
\begin{prop}\label{propdrift}
Let $p\in(2,\infty)$ and $\mathbf H\in L^2(U,\mathbb R^2)$. Then the drift term defines a bounded operator
\begin{equation} \label{driftoperat}
\mathcal{D}_{\mathbf H}:H^{1,p}_0(U)\longrightarrow H^{-1,p}(U),
\qquad
\langle\mathcal{D}_{\mathbf H}u,\varphi\rangle
:=\int_U \langle\mathbf H,\nabla u\rangle\varphi\,dx,
\end{equation}
for $u\in H^{1,p}_0(U)$ and $\varphi\in H^{1,p'}_0(U)$, where $p':=\frac{p}{p-1}$. More precisely, there exists a constant $C_D>0$, depending only on $p$, such that
\begin{equation}\label{driftbound}
\|\mathcal{D}_{\mathbf H}u\|_{H^{-1,p}(U)}
\le C_D\|\mathbf H\|_{L^2(U)}\|u\|_{H^{1,p}(U)}
\qquad\text{for all }u\in H^{1,p}_0(U).
\end{equation}
Moreover, for every $\eta>0$, choose $\mathbf H_\eta\in C_0^\infty(U,\mathbb R^2)$ such that
$$
C_D\|\mathbf H-\mathbf H_\eta\|_{L^2(U)}\le\eta,
$$
where $C_D>0$ is the constant in \eqref{driftbound}. Then
\begin{equation}\label{driftepsilon}
\begin{aligned}
\|\mathcal{D}_{\mathbf H}u\|_{H^{-1,p}(U)}
&\le \eta\|u\|_{H^{1,p}(U)}+
\bigl(
\|\mathbf H_\eta\|_{L^\infty(U)}
+\|\operatorname{div}\mathbf H_\eta\|_{L^\infty(U)}
\bigr)\|u\|_{L^p(U)}
\end{aligned}
\end{equation}
for all $u\in H^{1,p}_0(U)$.
\end{prop}
\begin{proof}
Set
$$
q:=\frac{2p}{p-2}.
$$
Since $\frac12+\frac1p+\frac1q=1$ and the two-dimensional Sobolev embedding gives
$$
H^{1,p'}_0(U)\hookrightarrow L^q(U),
$$
we have, for every $u\in H^{1,p}_0(U)$ and $\varphi\in H^{1,p'}_0(U)$,
\begin{align*}
\left|\int_U \langle\mathbf H,\nabla u\rangle\varphi\,dx\right|
&\le \|\mathbf H\|_{L^2(U)}
\|\nabla u\|_{L^p(U)}\|\varphi\|_{L^q(U)}\\
&\le C_D\|\mathbf H\|_{L^2(U)}
\|u\|_{H^{1,p}(U)}\|\varphi\|_{H^{1,p'}(U)}.
\end{align*}
This proves \eqref{driftbound}.\\
Let $\eta>0$. By the density of $C_0^\infty(U,\mathbb R^2)$ in $L^2(U,\mathbb R^2)$, choose $\mathbf H_\eta$ as in the statement. Applying \eqref{driftbound} to $\mathbf H-\mathbf H_\eta$, we obtain
$$
\|\mathcal{D}_{\mathbf H-\mathbf H_\eta}u\|_{H^{-1,p}(U)}
\le C_D\|\mathbf H-\mathbf H_\eta\|_{L^2(U)}
\|u\|_{H^{1,p}(U)}
\le \eta\|u\|_{H^{1,p}(U)}
\qquad\text{for all }u\in H^{1,p}_0(U).
$$
On the other hand, integration by parts gives
$$
\int_U \langle\mathbf H_\eta,\nabla u\rangle\varphi\,dx
=-\int_U u\langle\mathbf H_\eta,\nabla\varphi\rangle\,dx
-\int_U u(\operatorname{div}\mathbf H_\eta)\varphi\,dx
$$
for every $\varphi\in C_0^\infty(U)$. Hence
\begin{align*}
\left|\int_U \langle\mathbf H_\eta,\nabla u\rangle\varphi\,dx\right|
&\le \|\mathbf H_\eta\|_{L^\infty(U)}
\|u\|_{L^p(U)}\|\nabla\varphi\|_{L^{p'}(U)}\\
&\quad+
\|\operatorname{div}\mathbf H_\eta\|_{L^\infty(U)}
\|u\|_{L^p(U)}\|\varphi\|_{L^{p'}(U)}\\
&\le
\bigl(
\|\mathbf H_\eta\|_{L^\infty(U)}
+\|\operatorname{div}\mathbf H_\eta\|_{L^\infty(U)}
\bigr)
\|u\|_{L^p(U)}\|\varphi\|_{H^{1,p'}(U)}.
\end{align*}
By density, the same estimate holds for every $\varphi\in H^{1,p'}_0(U)$, and therefore
$$
\|\mathcal{D}_{\mathbf H_\eta}u\|_{H^{-1,p}(U)}
\le
\bigl(
\|\mathbf H_\eta\|_{L^\infty(U)}
+\|\operatorname{div}\mathbf H_\eta\|_{L^\infty(U)}
\bigr)\|u\|_{L^p(U)}
$$
for every $u\in H^{1,p}_0(U)$. Therefore
\begin{align*}
\|\mathcal{D}_{\mathbf H}u\|_{H^{-1,p}(U)}
&\le
\|\mathcal{D}_{\mathbf H-\mathbf H_\eta}u\|_{H^{-1,p}(U)}
+\|\mathcal{D}_{\mathbf H_\eta}u\|_{H^{-1,p}(U)}\\
&\le \eta\|u\|_{H^{1,p}(U)}+
\bigl(
\|\mathbf H_\eta\|_{L^\infty(U)}
+\|\operatorname{div}\mathbf H_\eta\|_{L^\infty(U)}
\bigr)\|u\|_{L^p(U)}
\end{align*}
for all $u\in H^{1,p}_0(U)$.
\end{proof}
\subsection{Uniform estimates along a homotopy}
The higher integrability of gradients for uniformly elliptic divergence-form equations with bounded measurable coefficients goes back to Meyers (\cite{Meyers63}). For the global Dirichlet problem on a bounded open set with Lipschitz boundary, Gr\"oger established a corresponding $W^{1,p}$-regularity theory in \cite{Groger89}. Related regularity results for exponents close to $2$ and elliptic equations with critical drifts have been obtained more recently (see, for instance, \cite{KwonJMAA}). We use the following consequence of Gr\"oger's result.
\begin{prop}[Gr\"oger]\label{propgroger}
Let $U\subset\mathbb R^2$ be a bounded open set with Lipschitz boundary and let $A$ satisfy \eqref{elliptbdda}. There exists $\varepsilon_0>0$, depending only on the ellipticity and boundedness constants of $A$ and the geometry of $U$, such that, for every
$$
2<p<2+\varepsilon_0,
$$
the operator
\begin{equation} \label{operatormatha}
\mathcal A_p:H^{1,p}_0(U)\longrightarrow H^{-1,p}(U),
\qquad
\mathcal A_pu:=-\operatorname{div}(A\nabla u),
\end{equation}
is a bounded isomorphism. For each such $p$, set
$$
C_M
:=
\|\mathcal A_p^{-1}\|_{\mathcal L(H^{-1,p}(U),H^{1,p}_0(U))}.
$$
Then $C_M<\infty$ and
\begin{equation}\label{grogerestimate}
\|u\|_{H^{1,p}(U)}
\le
C_M\|\mathcal A_pu\|_{H^{-1,p}(U)}
\qquad
\text{for all }u\in H^{1,p}_0(U).
\end{equation}
\end{prop}
\begin{proof}
Since $U$ is a bounded open set with Lipschitz boundary, it is regular in the sense of \cite[Definition~2]{Groger89}. Hence, by \cite[Theorem~3]{Groger89}, there exists $q>2$ such that $U\in R_q$. We apply Gr\"oger's result to the shifted operator
$$
\widetilde{\mathcal A}u:=-\operatorname{div}(A\nabla u)+u.
$$
By \cite[Remark~11]{Groger89}, the structural assumptions in \cite[(4.1)]{Groger89} are satisfied, since $A$ is bounded and uniformly elliptic and the zero-order coefficient is identically equal to $1$. By \cite[Lemma~1 and Theorem~1]{Groger89}, the contraction condition in that theorem holds for all $p\in(2,\infty)$ sufficiently close to $2$. Consequently, there exists $p_0\in(2,q]$ such that
$$
\widetilde{\mathcal A}:H^{1,p}_0(U)\longrightarrow H^{-1,p}(U)
$$
is a bounded isomorphism for every $p\in(2,p_0)$. Here, $p_0$ can be chosen depending only on $\lambda$, $M$, and the geometry of $U$. Set $\varepsilon_0:=p_0-2>0$.\\
We next remove the zero-order term from this isomorphism. Fix $p\in(2,2+\varepsilon_0)$ and $f\in H^{-1,p}(U)$. Since $U$ is bounded and $p'<2$, we have
$$
H^{-1,p}(U)\hookrightarrow H^{-1,2}(U).
$$
The Lax-Milgram theorem therefore yields a unique $u\in H^{1,2}_0(U)$ satisfying
$$
-\operatorname{div}(A\nabla u)=f.
$$
Moreover, the two-dimensional Sobolev embedding gives $H^{1,p'}_0(U)\hookrightarrow L^2(U)$. Thus $u\in L^2(U)$ defines an element of $H^{-1,p}(U)$. Consequently, $f+u\in H^{-1,p}(U)$. By the invertibility of $\widetilde{\mathcal A}$, there exists $v\in H^{1,p}_0(U)$ satisfying
$$
-\operatorname{div}(A\nabla v)+v=f+u.
$$
Since $H^{1,p}_0(U)\hookrightarrow H^{1,2}_0(U)$, both $u$ and $v$ solve the same shifted equation in $H^{1,2}_0(U)$. Indeed, their difference $w:=u-v\in H^{1,2}_0(U)$ satisfies
$$
-\operatorname{div}(A\nabla w)+w=0.
$$
Testing with $w$ and using ellipticity gives
$$
\lambda\int_U|\nabla w|^2\,dx+\int_U|w|^2\,dx
\le\int_U\langle A\nabla w,\nabla w\rangle\,dx
+\int_U|w|^2\,dx=0.
$$
Thus $w=0$, so $u=v\in H^{1,p}_0(U)$. This proves that $\mathcal A_p$ is surjective.
Its injectivity follows from the same energy argument.
The boundedness of $\mathcal A_p$ follows from $a_{ij}\in L^\infty(U)$ for all $1\le i,j\le2$ and H\"older's inequality. Therefore $\mathcal A_p$ is a bounded bijection from $H^{1,p}_0(U)$ onto $H^{-1,p}(U)$. By the bounded inverse theorem, $\mathcal A_p^{-1}$ is bounded, and \eqref{grogerestimate} follows.
\end{proof}
\\
Let $s\in(1,\infty)$ and $p\in(2,2+\varepsilon_0)$ satisfy
\begin{equation*}\label{padmissible}
s\ge\frac{2p}{p+2},
\end{equation*}
where $\varepsilon_0>0$ is as in Proposition~\ref{propgroger}. 
Let $\mathcal A_p$, $\mathcal D_{\mathbf{H}}$ and $\mathcal Z_c$ be the operators defined in \eqref{operatormatha}, \eqref{driftoperat} and \eqref{zeroordermult}, respectively. For each $t\in[0,1]$, define
\begin{equation}\label{Ttoperator}
\mathcal T_t:=\mathcal A_p+t\mathcal D_{\mathbf H}+t\mathcal Z_c,
\end{equation}
that is,
$$
\mathcal T_tu=-\operatorname{div}(A\nabla u)+t\langle\mathbf H,\nabla u\rangle+tcu, \quad u \in H^{1,p}_0(U).
$$
By Propositions~\ref{propdrift} and \ref{propzeroorder},
$$
\mathcal T_t:H^{1,p}_0(U)\longrightarrow H^{-1,p}(U)
$$
is bounded for every $t\in[0,1]$.
\begin{prop}[uniform a priori estimate]\label{uniformapriori}
Let $U\subset\mathbb R^2$ be a bounded open set with Lipschitz boundary and let $A=(a_{ij})_{1\le i,j\le 2}$ satisfy \eqref{elliptbdda}. Assume that
$$
\mathbf H\in L^2(U,\mathbb R^2),
\qquad
c\in L^s(U),
\qquad
c\ge0,
$$
for some $s\in(1,\infty)$. Let $\varepsilon_0>0$ be as in Proposition~\ref{propgroger} and $p\in(2,2+\varepsilon_0)$. Assume that
$$
s\ge\frac{2p}{p+2}.
$$
For $t\in[0,1]$, let $\mathcal T_t$ be defined by \eqref{Ttoperator}. Let
$$
C_M:=
\|\mathcal A_p^{-1}\|_{\mathcal L(H^{-1,p}(U),H^{1,p}_0(U))}
$$
and set
$$
\eta_0:=\frac{1}{4C_M}.
$$
Denote by $C_Z$ and $C_D$ the constants in \eqref{mcbounde} and \eqref{driftbound}, respectively. Choose $c_{\eta_0}\in L^\infty(U)$ and $\mathbf H_{\eta_0}\in C_0^\infty(U,\mathbb R^2)$ such that
$$
C_Z\|c-c_{\eta_0}\|_{L^s(U)}\le\eta_0,
\qquad
C_D\|\mathbf H-\mathbf H_{\eta_0}\|_{L^2(U)}\le\eta_0,
$$
and define
$$
C_{\eta_0}^{(1)}:=\|c_{\eta_0}\|_{L^\infty(U)},
\qquad
C_{\eta_0}^{(2)}:=
\|\mathbf H_{\eta_0}\|_{L^\infty(U)}
+\|\operatorname{div}\mathbf H_{\eta_0}\|_{L^\infty(U)}.
$$
Then
\begin{equation}\label{aprioriestimate}
\|u\|_{H^{1,p}(U)}
\le
2C_M\|\mathcal T_tu\|_{H^{-1,p}(U)}
+
2C_M\left(C_{\eta_0}^{(1)}+C_{\eta_0}^{(2)}\right)
\|u\|_{L^p(U)}
\end{equation}
for all $u\in H^{1,p}_0(U)$ and $t\in[0,1]$.
\end{prop}
\begin{proof}
By \eqref{grogerestimate} and the triangle inequality,
\begin{align*}
\|u\|_{H^{1,p}(U)}
&\le C_M\|\mathcal A_pu\|_{H^{-1,p}(U)}\\
&\le C_M\|\mathcal T_tu\|_{H^{-1,p}(U)}
+C_Mt\|\mathcal D_{\mathbf H}u\|_{H^{-1,p}(U)}
+C_Mt\|\mathcal Z_cu\|_{H^{-1,p}(U)}.
\end{align*}
By \eqref{driftepsilon} and \eqref{cepsilon}, with the choices of $\mathbf H_{\eta_0}$ and $c_{\eta_0}$ made above,
\begin{align*}
\|\mathcal D_{\mathbf H}u\|_{H^{-1,p}(U)}
&\le \eta_0\|u\|_{H^{1,p}(U)}
+C_{\eta_0}^{(2)}\|u\|_{L^p(U)},\\
\|\mathcal Z_cu\|_{H^{-1,p}(U)}
&\le \eta_0\|u\|_{H^{1,p}(U)}
+C_{\eta_0}^{(1)}\|u\|_{L^p(U)}.
\end{align*}
Since $0\le t\le1$ and $2C_M\eta_0=\frac12$, it follows that
$$
\|u\|_{H^{1,p}(U)}
\le
C_M\|\mathcal T_tu\|_{H^{-1,p}(U)}
+\frac12\|u\|_{H^{1,p}(U)}
+C_M\left(C_{\eta_0}^{(1)}+C_{\eta_0}^{(2)}\right)
\|u\|_{L^p(U)}.
$$
Subtracting $\frac12\|u\|_{H^{1,p}(U)}$ from both sides and multiplying by $2$ yields \eqref{aprioriestimate}.
\end{proof}
\begin{prop}[uniform inverse estimate]\label{uniforminverse}
Under the assumptions of Proposition~\ref{uniformapriori}, there exists a constant $C>0$, independent of $t\in[0,1]$ and $u \in H^{1,p}_0(U)$, such that
\begin{equation}\label{uniforminverseestimate}
\|u\|_{H^{1,p}(U)}
\le C\|\mathcal T_tu\|_{H^{-1,p}(U)}
\qquad\text{for all }u\in H^{1,p}_0(U)
\text{ and }t\in[0,1].
\end{equation}
\end{prop}
\begin{proof}
Suppose that \eqref{uniforminverseestimate} is false. Then, for every $n\ge1$, there exist $t_n\in[0,1]$ and $w_n\in H^{1,p}_0(U)$ such that
$$
\|w_n\|_{H^{1,p}(U)}
>n\|\mathcal T_{t_n}w_n\|_{H^{-1,p}(U)}.
$$
In particular, $w_n\ne0$. Setting
$$
u_n:=\frac{w_n}{\|w_n\|_{H^{1,p}(U)}},
$$
we obtain
\begin{equation}\label{normalized}
\|u_n\|_{H^{1,p}(U)}=1,
\qquad
\|\mathcal T_{t_n}u_n\|_{H^{-1,p}(U)}<\frac1n.
\end{equation}
Since $[0,1]$ is compact, after passing to a subsequence we may assume that
$$
\lim_{n\to\infty}t_n=t
\qquad\text{for some }t\in[0,1].
$$
Moreover, since $H^{1,p}_0(U)$ is reflexive and $(u_n)_{n\ge1}$ is bounded in $H^{1,p}_0(U)$, after passing to a further subsequence there exists $u\in H^{1,p}_0(U)$ such that
\begin{equation}\label{weakun}
\lim_{n\to\infty}u_n=u
\qquad\text{weakly in }H^{1,p}_0(U).
\end{equation}
By the Rellich-Kondrachov theorem, after passing to a further subsequence if necessary, we also have
\begin{equation}\label{strongun}
\lim_{n\to\infty}u_n=u
\qquad\text{in }L^p(U).
\end{equation}
The limit in \eqref{strongun} agrees with the weak limit in \eqref{weakun}, since the embedding of $H^{1,p}_0(U)$ into $L^p(U)$ is continuous.\\
We first show that $u\ne0$. Indeed, if $u=0$, then \eqref{strongun} gives
$$
\lim_{n\to\infty}\|u_n\|_{L^p(U)}=0.
$$
Applying Proposition~\ref{uniformapriori} to $u_n$ and $t_n$, we obtain
$$
\begin{aligned}
1=\|u_n\|_{H^{1,p}(U)}
&\le
2C_M\|\mathcal T_{t_n}u_n\|_{H^{-1,p}(U)}+
2C_M\left(C_{\eta_0}^{(1)}+C_{\eta_0}^{(2)}\right)
\|u_n\|_{L^p(U)}.
\end{aligned}
$$
By \eqref{normalized} and \eqref{strongun}, the right-hand side has limit zero as $n$ tends to infinity, which is impossible. Hence
\begin{equation}\label{nonzerolimit}
u\ne0.
\end{equation}
On the other hand, since $\mathcal A_p$, $\mathcal D_{\mathbf H}$, and $\mathcal Z_c$ are bounded linear operators, \eqref{weakun} implies
$$
\lim_{n\to\infty}\mathcal A_pu_n=\mathcal A_pu,
\qquad
\lim_{n\to\infty}\mathcal D_{\mathbf H}u_n
=\mathcal D_{\mathbf H}u,
\qquad
\lim_{n\to\infty}\mathcal Z_cu_n=\mathcal Z_cu
$$
weakly in $H^{-1,p}(U)$. In particular,
$$
\lim_{n\to\infty}\mathcal T_tu_n=\mathcal T_tu
\qquad\text{weakly in }H^{-1,p}(U).
$$
To account for the varying parameter $t_n$, we use
$$
\|\mathcal T_{t_n}u_n-\mathcal T_tu_n\|_{H^{-1,p}(U)}
\le
|t_n-t|
\left(
\|\mathcal D_{\mathbf H}u_n\|_{H^{-1,p}(U)}
+\|\mathcal Z_cu_n\|_{H^{-1,p}(U)}
\right).
$$
Since $\lim_{n\to\infty}t_n=t$, $(u_n)_{n\ge1}$ is bounded in $H^{1,p}_0(U)$, and $\mathcal D_{\mathbf H}$ and $\mathcal Z_c$ are bounded operators, it follows that
$$
\lim_{n\to\infty}
\|\mathcal T_{t_n}u_n-\mathcal T_tu_n\|_{H^{-1,p}(U)}
=0.
$$
Therefore,
$$
\lim_{n\to\infty}\mathcal T_{t_n}u_n=\mathcal T_tu
\qquad\text{weakly in }H^{-1,p}(U).
$$
But \eqref{normalized} gives
$$
\lim_{n\to\infty}\mathcal T_{t_n}u_n=0
\qquad\text{in }H^{-1,p}(U).
$$
Strong convergence implies weak convergence to the same limit. By uniqueness of the weak limit, we conclude that
$$
\mathcal T_tu=0.
$$
Since
$$
t\mathbf H\in L^2(U,\mathbb R^2),
\qquad
tc\in L^s(U)\subset L^1(U),
\qquad
tc\ge0,
$$
Corollary~\ref{coruniqueness} yields $u=0$, which contradicts \eqref{nonzerolimit}. Hence \eqref{uniforminverseestimate} holds.
\end{proof}
\begin{lem}\label{smalloperator}
Let $X$ be a Banach space and let
$$
T:X\longrightarrow X
$$
be a bounded linear operator satisfying
$$
\|T\|_{\mathcal L(X)}<1.
$$
Then
$$
I+T:X\longrightarrow X
$$
is a bounded isomorphism.
\end{lem}
\begin{proof}
If $(I+T)x=0$, then
$$
\|x\|=\|Tx\|\le\|T\|_{\mathcal L(X)}\|x\|,
$$
and hence $x=0$. Thus $I+T$ is injective.\\
For $y\in X$, consider
$$
\Phi(x):=y-Tx.
$$
Since
$$
\|\Phi(x)-\Phi(z)\|
\le
\|T\|_{\mathcal L(X)}\|x-z\|,
$$
$\Phi$ is a contraction on $X$. By the Banach fixed point theorem, there exists a unique $x\in X$ such that
$$
x=\Phi(x)=y-Tx,
$$
or equivalently,
$$
(I+T)x=y.
$$
Thus $I+T$ is bijective. Moreover,
$$
\|(I+T)x\|
\ge
\left(1-\|T\|_{\mathcal L(X)}\right)\|x\|,
$$
so
$$
\|(I+T)^{-1}\|_{\mathcal L(X)}
\le
\frac{1}{1-\|T\|_{\mathcal L(X)}}.
$$
\end{proof}
\\
The proof of the following theorem combines Gr\"oger's isomorphism for $p\in(2,\infty)$ sufficiently close to $2$ with the lower-order estimates and the uniqueness result established above. Related arguments based on $\varepsilon$-estimates, compactness, and uniqueness appear in \cite{Kwon}. The use of lower-order estimates together with the method of continuity can also be found in \cite{KimOh23}. In the present setting, compactness and Corollary~\ref{coruniqueness} yield a uniform estimate along the family $\mathcal T_t$, which allows us to apply the method of continuity.
\begin{theo}[$H^{1,p}$ well-posedness]\label{thmexistence}
Under the assumptions of Proposition~\ref{uniformapriori}, the operator
$$
\mathcal T_1
=\mathcal A_p+\mathcal D_{\mathbf H}+\mathcal Z_c:
H^{1,p}_0(U)\longrightarrow H^{-1,p}(U)
$$
is a bounded isomorphism. Consequently, for every $\ell\in H^{-1,p}(U)$, there exists a unique $u\in H^{1,p}_0(U)$ satisfying
\begin{equation}\label{generalrhs}
-\operatorname{div}(A\nabla u)
+\langle\mathbf H,\nabla u\rangle+cu
=\ell
\qquad\text{in }H^{-1,p}(U),
\end{equation}
and
\begin{equation}\label{generalest}
\|u\|_{H^{1,p}(U)}
\le C\|\ell\|_{H^{-1,p}(U)},
\end{equation}
where $C>0$ is independent of $\ell$. In particular, if
$$
f\in L^r(U),
\qquad
\mathbf F\in L^p(U,\mathbb R^2),
\qquad
\frac{2p}{p+2}\le r\le\infty,
$$
then
$$
\ell:=f-\operatorname{div}\mathbf F\in H^{-1,p}(U),
$$
and the unique solution $u$ satisfies
\begin{equation}\label{frest}
\|u\|_{H^{1,p}(U)}
\le C\left(
\|f\|_{L^r(U)}+\|\mathbf F\|_{L^p(U)}
\right),
\end{equation}
with a constant $C>0$ independent of $f$ and $\mathbf F$.\\
For each fixed $s\in(1,\infty)$, all the preceding conclusions hold for every
$$
2<p<2+\varepsilon_s,
$$
where
$$
\varepsilon_s:=
\begin{cases}
\displaystyle
\min\left\{
\varepsilon_0,\frac{4(s-1)}{2-s}
\right\},
& s\in(1,2),\\[1.2ex]
\varepsilon_0,
& s\in[2,\infty).
\end{cases}
$$
\end{theo}
\begin{proof}
For $t\in[0,1]$, set
$$
\mathcal T_t
:=\mathcal A_p+t\mathcal D_{\mathbf H}+t\mathcal Z_c.
$$
By Propositions~\ref{propgroger}, \ref{propdrift}, and \ref{propzeroorder}, each operator
$$
\mathcal T_t:H^{1,p}_0(U)\longrightarrow H^{-1,p}(U)
$$
is bounded. Moreover, Corollary~\ref{coruniqueness} implies that $\mathcal T_t$ is injective for every $t\in[0,1]$. Define
$$
\mathcal S:=
\left\{
t\in[0,1]:\mathcal T_t\text{ is surjective}
\right\}.
$$
We claim that $\mathcal S=[0,1]$. By Proposition~\ref{propgroger},
$$
\mathcal T_0=\mathcal A_p
$$
is an isomorphism. Hence $0\in\mathcal S$.\\
We next show that $\mathcal S$ is relatively open in $[0,1]$. Fix $t_0\in\mathcal S$ and set
$$
\mathcal P:=\mathcal D_{\mathbf H}+\mathcal Z_c.
$$
Since $\mathcal T_{t_0}$ is a bounded bijection, the bounded inverse theorem implies that $\mathcal T_{t_0}^{-1}$ is bounded. Moreover,
$$
\mathcal T_t
=\mathcal T_{t_0}
\left[
I+(t-t_0)\mathcal T_{t_0}^{-1}\mathcal P
\right],
$$
where $I$ denotes the identity operator on $H^{1,p}_0(U)$. If $t\in[0,1]$ satisfies
$$
|t-t_0|
<
\frac{1}{
\|\mathcal T_{t_0}^{-1}\mathcal P\|_{\mathcal L(H^{1,p}_0(U))}
+1},
$$
then
$$
\left\|
(t-t_0)\mathcal T_{t_0}^{-1}\mathcal P
\right\|_{\mathcal L(H^{1,p}_0(U))}
<1.
$$
By Lemma~\ref{smalloperator}, the operator
$$
I+(t-t_0)\mathcal T_{t_0}^{-1}\mathcal P
$$
is invertible. Consequently, $\mathcal T_t$ is invertible and, in particular, surjective. Thus $\mathcal S$ is relatively open in $[0,1]$.\\
To show that $\mathcal S$ is relatively closed in $[0,1]$, let $(t_n)_{n\ge1}$ be a sequence in $\mathcal S$ such that
$$
\lim_{n\to\infty}t_n=t
\qquad\text{for some }t\in[0,1],
$$
and fix $\ell\in H^{-1,p}(U)$. For each $n\ge1$, there exists $u_n\in H^{1,p}_0(U)$ such that
$$
\mathcal T_{t_n}u_n=\ell.
$$
By Proposition~\ref{uniforminverse},
$$
\|u_n\|_{H^{1,p}(U)}
\le C\|\ell\|_{H^{-1,p}(U)},
$$
where $C$ is independent of $n$ and $\ell$. Since $H^{1,p}_0(U)$ is reflexive, after passing to a subsequence we may assume that for some $u\in H^{1,p}_0(U)$,
$$
\lim_{n\to\infty}u_n=u
\qquad\text{weakly in }H^{1,p}_0(U).
$$
 Since $\mathcal A_p$, $\mathcal D_{\mathbf H}$, and $\mathcal Z_c$ are bounded linear operators,
$$
\lim_{n\to\infty}\mathcal A_pu_n=\mathcal A_pu,
\qquad
\lim_{n\to\infty}\mathcal D_{\mathbf H}u_n
=\mathcal D_{\mathbf H}u,
\qquad
\lim_{n\to\infty}\mathcal Z_cu_n=\mathcal Z_cu
$$
weakly in $H^{-1,p}(U)$. In particular,
$$
\lim_{n\to\infty}\mathcal T_tu_n=\mathcal T_tu
\qquad\text{weakly in }H^{-1,p}(U).
$$
Moreover, since $\lim_{n\to\infty}t_n=t$ and $(\mathcal Pu_n)_{n\ge1}$ is bounded in $H^{-1,p}(U)$,
$$
\begin{aligned}
\lim_{n\to\infty}
\|\mathcal T_{t_n}u_n-\mathcal T_tu_n\|_{H^{-1,p}(U)}
&=
\lim_{n\to\infty}
|t_n-t|\|\mathcal Pu_n\|_{H^{-1,p}(U)}\\
&=0.
\end{aligned}
$$
It follows that
$$
\lim_{n\to\infty}\mathcal T_{t_n}u_n=\mathcal T_tu
\qquad\text{weakly in }H^{-1,p}(U).
$$
Since $\mathcal T_{t_n}u_n=\ell$ for every $n\ge1$, we obtain
$$
\mathcal T_tu=\ell.
$$
As $\ell\in H^{-1,p}(U)$ was arbitrary, $\mathcal T_t$ is surjective and hence $t\in\mathcal S$. Thus $\mathcal S$ is relatively closed in $[0,1]$.\\
Since $[0,1]$ is connected and $\mathcal S$ is nonempty, relatively open, and relatively closed in $[0,1]$, we conclude that
$$
\mathcal S=[0,1].
$$
In particular, $\mathcal T_1$ is surjective. Together with its injectivity, this shows that $\mathcal T_1$ is a bounded bijection. By the bounded inverse theorem, $\mathcal T_1$ is a bounded isomorphism, and \eqref{generalest} follows.\\
It remains to verify the assertion for right-hand sides of the form $f-\operatorname{div}\mathbf F$. Let $r'$ denote the H\"older conjugate of $r$, with $r'=1$ when $r=\infty$. The condition
$$
r\ge\frac{2p}{p+2}
$$
implies that
$$
r'\le\frac{2p}{p-2}.
$$
Since $1<p'<2$ and $U$ is bounded, the two-dimensional Sobolev embedding gives
$$
H^{1,p'}_0(U)
\hookrightarrow
L^{\frac{2p}{p-2}}(U)
\hookrightarrow
L^{r'}(U).
$$
Hence, for every $\varphi\in H^{1,p'}_0(U)$,
\begin{align*}
\left|
\int_U f\varphi\,dx
+\int_U\langle\mathbf F,\nabla\varphi\rangle\,dx
\right|
&\le
\|f\|_{L^r(U)}\|\varphi\|_{L^{r'}(U)}
+\|\mathbf F\|_{L^p(U)}
\|\nabla\varphi\|_{L^{p'}(U)}\\
&\le
C_1\left(
\|f\|_{L^r(U)}+\|\mathbf F\|_{L^p(U)}
\right)\|\varphi\|_{H^{1,p'}(U)},
\end{align*}
where $C_1>0$ depends only on $p,r,|U|$. Therefore,
$$
f-\operatorname{div}\mathbf F\in H^{-1,p}(U)
$$
and
\begin{equation}\label{fboldest}
\|f-\operatorname{div}\mathbf F\|_{H^{-1,p}(U)}
\le
C_1\left(
\|f\|_{L^r(U)}+\|\mathbf F\|_{L^p(U)}
\right).
\end{equation}
Combining \eqref{fboldest} with \eqref{generalest} gives \eqref{frest}.\\
Finally, if $s\in(1,2)$, then
$$
s\ge\frac{2p}{p+2}
$$
is equivalent to
$$
p\le\frac{2s}{2-s}
=2+\frac{4(s-1)}{2-s}.
$$
Thus the stated choice of $\varepsilon_s$ guarantees both $p<2+\varepsilon_0$ and $s\ge\frac{2p}{p+2}$ whenever $2<p<2+\varepsilon_s$. If $s\ge2$, the condition $s\ge\frac{2p}{p+2}$ holds automatically for every $p\in(2,\infty)$, so the only restriction is $2<p<2+\varepsilon_0$. This completes the proof.
\end{proof}
\\
\begin{proof}{\bf of Theorem~\ref{thmmain}~{\rm(ii)}.}
Let $\varepsilon_s>0$ be as in Theorem~\ref{thmexistence}. Fix $q\in(2,\infty)$ and set
$$
\varepsilon:=\frac12\min\{\varepsilon_s,q-2\},
\qquad
p:=2+\varepsilon.
$$
Then
$$
0<\varepsilon<q-2,
\qquad
2<p<2+\varepsilon_s,
\qquad
p<q.
$$
By the definition of $\varepsilon_s$ and Proposition~\ref{propgroger}, $\varepsilon$ depends only on $q$, $s$, the constants $\lambda,M$ in \eqref{elliptbdda}, and the geometry of $U$.\\
Let $\ell\in H^{-1,q}(U)$. Since $U$ is bounded and $q'<p'$, H\"older's inequality gives the continuous embedding
$$
H^{1,p'}_0(U)\hookrightarrow H^{1,q'}_0(U).
$$
Indeed, the corresponding norm estimate holds for smooth functions, and approximation by functions in $C_0^\infty(U)$ preserves the zero boundary condition. Consequently, for every $\varphi\in H^{1,p'}_0(U)$,
$$
|\langle\ell,\varphi\rangle|
\le
\|\ell\|_{H^{-1,q}(U)}
\|\varphi\|_{H^{1,q'}(U)}
\le
C_1\|\ell\|_{H^{-1,q}(U)}
\|\varphi\|_{H^{1,p'}(U)},
$$
where $C_1>0$ depends only on $p$, $q$, and $|U|$. Thus the restriction of $\ell$ to $H^{1,p'}_0(U)$ defines an element of $H^{-1,p}(U)$, and
$$
\|\ell\|_{H^{-1,p}(U)}
\le C_1\|\ell\|_{H^{-1,q}(U)}.
$$
We identify $\ell$ with this restriction, which agrees with the original functional on $C_0^\infty(U)$.\\
Since $2<p<2+\varepsilon_s$, Theorem~\ref{thmexistence} yields a unique $u\in H^{1,p}_0(U)$ satisfying \eqref{generalrhs}. Testing this identity against $\varphi\in C_0^\infty(U)$ gives \eqref{introweakformulation}, so $u$ is an $H^{1,p}$-weak solution of \eqref{introproblem}. Furthermore, \eqref{generalest} gives
$$
\|u\|_{H^{1,p}(U)}
\le C_2\|\ell\|_{H^{-1,p}(U)}
\le C_1C_2\|\ell\|_{H^{-1,q}(U)},
$$
where $C_2>0$ is independent of $\ell$. This proves \eqref{mainestimate}. Uniqueness in the class of $H^{1,p}$-weak solutions follows from part {\rm(i)} by linearity.
\end{proof}
\subsection{Failure of $H^{1,p}$ solvability at the endpoint $c\in L^1$}
We now show that the assumption $s>1$ in Theorem~\ref{thmexistence} cannot in general be replaced by $s=1$, even when the principal part is the Laplacian and the right-hand side is smooth.
\begin{prop}[endpoint counterexample]\label{endpointcounterexample}
Let $B_1:=B_1(0)\subset\mathbb R^2$. There exist a coefficient
$$
c\in L^1(B_1),
\qquad
c\ge0,
$$
and a function $f\in C_0^\infty(B_1)$ such that
\begin{equation}\label{counterPDE}
-\Delta u+cu=f
\qquad\text{in }B_1
\end{equation}
admits a solution $u\in H^{1,2}_0(B_1)\cap L^\infty(B_1)$, but has no solution in $H^{1,p}_0(B_1)$ for any $p\in(2,\infty)$. Here the equation is understood in the sense of distributions.
\end{prop}
\begin{proof}
Fix $\beta>0$ and set
$$
\tau(r):=\ln\frac er,
\qquad 0<r<1.
$$
Define
\begin{equation*}\label{counterc}
c(x):=
\frac{\beta(\beta+1)}
{\|x\|^2\left(\ln\frac e{\|x\|}\right)^2}
\qquad\text{for }x\in B_1\setminus\{0\},
\end{equation*}
and set $c(0):=0$. Then $c\ge0$. Using polar coordinates and the substitution $t=\ln(\frac{e}{r})$, for which $r=e^{1-t}$ and $\frac{1}{r}dr=-dt$, we obtain
$$
\begin{aligned}
\int_{B_1}c(x)\,dx
&=
2\pi\beta(\beta+1)
\int_0^1\frac{dr}{r(\ln(\frac{e}{r}))^2}=
2\pi\beta(\beta+1)
\int_1^\infty\frac{dt}{t^2}
=
2\pi\beta(\beta+1)
<\infty.
\end{aligned}
$$
Hence $c\in L^1(B_1)$. For every $s>1$, a similar calculation gives
$$
\begin{aligned}
\int_{B_1}c(x)^s\,dx
&=
2\pi[\beta(\beta+1)]^s
\int_0^1r^{1-2s}\left(\ln\frac er\right)^{-2s}\,dr=
2\pi[\beta(\beta+1)]^s e^{2-2s}
\int_1^\infty e^{(2s-2)t}t^{-2s}\,dt=\infty.
\end{aligned}
$$
Thus $c\notin L^s(B_1)$ for every $s>1$.\\
Choose a radial function $\eta\in C_0^\infty(B_1)$ such that $0\le\eta\le1$ and
$$
\eta\equiv1\quad\text{on }B_{\frac{1}{4}}, \qquad
\eta\equiv0\quad\text{on }B_1\setminus B_{\frac{1}{2}}.
$$
Define
\begin{equation*}\label{counteru}
u(x):=
\eta(x)\left(\ln\frac e{\|x\|}\right)^{-\beta}
\qquad\text{for }x\in B_1\setminus\{0\},
\end{equation*}
and set $u(0):=0$. Since $\tau(r)>1$ for $0<r<1$, we have $0\le u\le1$. Moreover,
$$
\lim_{r\rightarrow0+}\tau(r)^{-\beta}=0.
$$
Thus $u$ is bounded and continuous on $B_1$, and its support is contained in $\overline{B}_{\frac{1}{2}}$.\\
For $0<r<1/4$, the radial profile of $u$ is $u(r)=\tau(r)^{-\beta}$. Since $\tau'(r)=-1/r$, differentiation gives
$$
u_r(r)
=-\beta\tau(r)^{-\beta-1}\tau'(r)
=\frac{\beta}{r\tau(r)^{\beta+1}}
$$
and
$$
u_{rr}(r)
=
-\frac{\beta}{r^2\tau(r)^{\beta+1}}
+
\frac{\beta(\beta+1)}{r^2\tau(r)^{\beta+2}}.
$$
Consequently,
$$
\Delta u
=
u_{rr}+\frac1r u_r
=
\frac{\beta(\beta+1)}{r^2\tau(r)^{\beta+2}}
=
c(r)u(r),
$$
and hence
\begin{equation}\label{countercancel}
-\Delta u+cu=0
\qquad\text{in }B_{\frac{1}{4}}\setminus\{0\}.
\end{equation}
We first verify that $u\in H^{1,2}_0(B_1)$. For $0<\|x\|<\frac{1}{4}$, the gradient of a radial function is $u_r(\|x\|)\frac{x}{\|x\|}$. Thus
$$
\nabla u(x)
=
\frac{\beta x}{\|x\|^2\left(\ln(\frac{e}{\|x\|})\right)^{\beta+1}}.
$$
It follows that
$$
\begin{aligned}
\int_{B_{1/4}}\|\nabla u\|^2\,dx
&=2\pi\beta^2
\int_0^{1/4}
\frac{dr}{r(\ln(\frac{e}{r}))^{2\beta+2}}\\
&=
2\pi\beta^2
\int_{\ln(4e)}^\infty t^{-2\beta-2}\,dt\\
&=
\frac{2\pi\beta^2}{2\beta+1}
\bigl(\ln(4e)\bigr)^{-2\beta-1}
<\infty.
\end{aligned}
$$
Away from the origin, $u$ is smooth and vanishes outside $B_{1/2}$. Its classical gradient on $B_1\setminus\{0\}$ therefore belongs to $L^2(B_1,\mathbb R^2)$.\\
To see that this is also the weak gradient on the whole ball, let $\varphi\in C_0^\infty(B_1)$ and $i\in\{1,2\}$. Integration by parts over $B_1\setminus\overline{B}_\rho$ gives
$$
\int_{B_1\setminus\overline{B}_\rho}
 u\,\partial_i\varphi\,dx
=
-\int_{B_1\setminus\overline{B}_\rho}
 (\partial_i u)\varphi\,dx
-\int_{\partial B_\rho}
 u\varphi\frac{x_i}{\rho}\,d\sigma.
$$
Moreover,
$$
\left|
\int_{\partial B_\rho}u\varphi\frac{x_i}{\rho}\,d\sigma
\right|
\le
2\pi\rho\|u\|_{L^\infty(B_1)}\|\varphi\|_{L^\infty(B_1)},
$$
and the limit of this upper bound is zero as $\rho$ decreases to zero. Both volume integrands are integrable, so passing to the limit yields
$$
\int_{B_1}u\,\partial_i\varphi\,dx
=
-\int_{B_1}(\partial_i u)\varphi\,dx.
$$
Since $u$ is bounded and $B_1$ has finite measure, $u\in L^2(B_1)$. Together with the gradient estimate, this proves that $u\in H^{1,2}(B_1)$. In addition, $u$ vanishes in a neighborhood of $\partial B_1$. Its extension by zero to $\mathbb R^2$ belongs to $H^{1,2}(\mathbb R^2)$, and convolution with smooth mollifiers of sufficiently small support gives functions in $C_0^\infty(B_1)$ converging to $u$ in $H^{1,2}(B_1)$. By the definition of $H^{1,2}_0(B_1)$, we conclude that
$$
u\in H^{1,2}_0(B_1).
$$
Moreover, since $c\in L^1(B_1)$ and $u$ is bounded,
$$
cu\in L^1(B_1),
\qquad
cu^2\in L^1(B_1).
$$
We next verify that \eqref{countercancel} holds distributionally across the origin. Let $\varphi\in C_0^\infty(B_{1/4})$. Multiplying \eqref{countercancel} by $\varphi$ and integrating by parts over $B_{1/4}\setminus\overline{B_\rho}$ gives
$$
\int_{B_{1/4}\setminus\overline{B}_\rho}
\langle\nabla u,\nabla\varphi\rangle\,dx
+
\int_{B_{1/4}\setminus\overline{B}_\rho}
cu\varphi\,dx
=
-u_r(\rho)\int_{\partial B_\rho}\varphi\,d\sigma.
$$
Here the boundary term on $\partial B_{1/4}$ is zero because $\varphi$ has compact support in $B_{1/4}$. On the inner boundary, the normal derivative of $u$ is $-u_r(\rho)$. Hence
$$
\left|
 u_r(\rho)\int_{\partial B_\rho}\varphi\,d\sigma
\right|
\le
2\pi\rho|u_r(\rho)|\|\varphi\|_{L^\infty(B_{1/4})}.
$$
The formula for $u_r$ shows that
$$
\lim_{\rho\rightarrow0+}\rho u_r(\rho)
=
\lim_{\rho\rightarrow0+}
\beta\tau(\rho)^{-\beta-1}
=0.
$$
The first volume integral has a limit because $\nabla u\in L^2(B_1,\mathbb R^2)$ and $\nabla\varphi$ is bounded. The second has a limit because $cu\in L^1(B_1)$. We therefore obtain
$$
\int_{B_{1/4}}\langle\nabla u,\nabla\varphi\rangle\,dx
+
\int_{B_{1/4}}cu\varphi\,dx
=0.
$$
Thus no additional distribution supported at the origin occurs, and in particular no Dirac mass is produced there. We have proved that
$$
-\Delta u+cu=0
\qquad\text{in }\mathcal D'(B_{1/4}).
$$
Define
$$
f:=-\Delta u+cu
$$
in the sense of distributions on $B_1$. By the preceding argument, $f$ vanishes on $B_{1/4}$. To identify it elsewhere, set
$$
g(x):=\left(\ln\frac e{\|x\|}\right)^{-\beta}
\qquad\text{for }x\in B_1\setminus\{0\}.
$$
The same radial computation as above gives $-\Delta g+cg=0$ on $B_1\setminus\{0\}$. Since $u=\eta g$, the product rule yields
$$
-\Delta u+cu
=
\eta(-\Delta g+cg)
-2\langle\nabla\eta,\nabla g\rangle
-g\Delta\eta
=
-2\langle\nabla\eta,\nabla g\rangle-g\Delta\eta
$$
on $B_1\setminus\{0\}$. The derivatives of $\eta$ vanish on $B_{1/4}$ and outside $\overline{B}_{1/2}$, while $g$ is smooth away from the origin. Thus the last expression extends by zero to a smooth function on $B_1$, and it represents the distribution $f$. In particular,
$$
\operatorname{supp}f
\subset
\left\{x\in B_1:\frac14\le\|x\|\le\frac12\right\},
\qquad
f\in C_0^\infty(B_1).
$$
We now check that the constructed solution has no higher Sobolev regularity. For every $p\in(2,\infty)$, the formula for the gradient gives
$$
\int_{B_{1/4}}\|\nabla u\|^p\,dx
=
2\pi\beta^p
\int_0^{1/4}
r^{1-p}\left(\ln\frac er\right)^{-p(\beta+1)}\,dr.
$$
With $t=\ln(\frac{e}{r})$, we obtain
$$
\int_{B_{1/4}}\|\nabla u\|^p\,dx
=
2\pi\beta^p e^{2-p}
\int_{\ln(4e)}^\infty
e^{(p-2)t}t^{-p(\beta+1)}\,dt.
$$
Since $p-2>0$ and
$$
\lim_{t\to\infty}\frac{\ln t}{t}=0,
$$
there exists $T>\ln(4e)$ such that
$$
p(\beta+1)\ln t\le\frac{p-2}{2}t
\qquad\text{for every }t\ge T.
$$
Consequently,
$$
e^{(p-2)t}t^{-p(\beta+1)}
\ge e^{(p-2)t/2}
\qquad\text{for every }t\ge T,
$$
so the integral diverges. Hence
\begin{equation}\label{counternotHp}
u\notin H^{1,p}(B_1)
\qquad\text{for every }p\in(2,\infty).
\end{equation}
It remains to exclude the existence of another solution with higher Sobolev regularity. Suppose that, for some $p\in(2,\infty)$, there exists $v\in H^{1,p}_0(B_1)$ satisfying
$$
-\Delta v+cv=f
\qquad\text{in }\mathcal D'(B_1).
$$
Since $p>2$, the Sobolev embedding gives $v\in L^\infty(B_1)$. Also, $H^{1,p}_0(B_1)\subset H^{1,2}_0(B_1)$ because $B_1$ has finite measure. Therefore
$$
w:=v-u\in H^{1,2}_0(B_1)\cap L^\infty(B_1)
$$
and
\begin{equation}\label{counterhomogeneous}
-\Delta w+cw=0
\qquad\text{in }\mathcal D'(B_1).
\end{equation}
Since $c\in L^1(B_1)$ and $w$ is bounded,
$$
cw\in L^1(B_1),
\qquad
cw^2\in L^1(B_1).
$$
By \cite[Proposition~A.8]{L25JMAA}, there exist a constant $M>0$ and a sequence $(\varphi_n)_{n\ge1}$ in $C_0^\infty(B_1)$ such that
$$
\sup_{n\ge1}\|\varphi_n\|_{L^\infty(B_1)}\le M,
\qquad
\lim_{n\to\infty}\varphi_n=w
\quad\text{in }H^{1,2}_0(B_1),
$$
and
$$
\lim_{n\to\infty}\varphi_n=w
\qquad\text{a.e. in }B_1.
$$
Testing \eqref{counterhomogeneous} with $\varphi_n$ gives
$$
\int_{B_1}\langle\nabla w,\nabla\varphi_n\rangle\,dx
+\int_{B_1}cw\varphi_n\,dx=0.
$$
The convergence in $H^{1,2}_0(B_1)$ implies
$$
\lim_{n\to\infty}\int_{B_1}\langle\nabla w,\nabla\varphi_n\rangle\,dx
=\int_{B_1}\|\nabla w\|^2\,dx.
$$
Moreover,
$$
|cw\varphi_n|\le M c|w|
\qquad\text{a.e. in }B_1,
$$
and $c|w|\in L^1(B_1)$. Hence the dominated convergence theorem gives
$$
\lim_{n\to\infty}\int_{B_1}cw\varphi_n\,dx
=\int_{B_1}cw^2\,dx.
$$
Passing to the limit in the tested equation, we obtain
$$
\int_{B_1}\|\nabla w\|^2\,dx+\int_{B_1}cw^2\,dx=0.
$$
Since $c\ge0$, both terms are nonnegative. Thus $\nabla w=0$ a.e. in $B_1$, and the Poincar\'e inequality yields $w=0$ because $w\in H^{1,2}_0(B_1)$. Therefore $v=u$, which contradicts \eqref{counternotHp}. Hence \eqref{counterPDE} has no solution in $H^{1,p}_0(B_1)$ for any $p\in(2,\infty)$.
\end{proof}

\section*{Acknowledgments}
This work was supported by the National Research Foundation of Korea (NRF) grant funded by the Korea government (MSIT) (RS-2026-25598810) and the Gyeongbuk ANCHOR system-(Regional Growth Innovation LAB) program through the Gyeongbuk ANCHOR Center, funded by the Ministry of Education (MOE) and the Gyeongsangbuk-do, Republic of Korea (2026-anchor-15-105).

\centerline{}
\centerline{}
Haesung Lee \\
Department of Mathematics and Big Data Science,  \\
Kumoh National Institute of Technology, \\
61 Daehak-ro, Gumi, Gyeongsangbuk-do 39177, Republic of Korea, \\
E-mail: fthslt@kumoh.ac.kr, \; fthslt14@gmail.com
\end{document}